\documentclass[11pt,reqno]{amsart}

\usepackage{amssymb,amsmath,epsfig}
\usepackage{amsfonts}
\usepackage{pst-plot}

\newtheorem{theorem}{Theorem}[section]

\newtheorem{lemma}[theorem]{Lemma}
\newtheorem{corollary}[theorem]{Corollary}
\newtheorem{proposition}[theorem]{Proposition}
\newtheorem{definition}[theorem]{Definition}
\newtheorem*{theorem*}{Theorem}

\numberwithin{equation}{section}

\title{Overpartition pairs and two classes of basic hypergeometric series}

\newcommand{\qbinom}[2]{\begin{bmatrix} #1 \\ #2 \end{bmatrix}_q}

\date{\today}
\author{Jeremy Lovejoy and Olivier Mallet}
\address{CNRS, LIAFA, Universit\'e Denis Diderot,
2, Place Jussieu, Case 7014, F-75251 Paris Cedex 05, FRANCE}
\email{lovejoy@liafa.jussieu.fr}
\address{LIAFA, Universit\'e Denis Diderot,
2, Place Jussieu, Case 7014, F-75251 Paris Cedex 05, FRANCE}
\email{mallet@liafa.jussieu.fr}


\begin{document}

\begin{abstract}
We study the combinatorics of two classes of basic hypergeometric
series.  We first show that these series are the generating
functions for certain overpartition pairs defined by frequency
conditions on the parts.  We then show that when specialized these
series are also the generating functions for overpartition pairs
with bounded successive ranks, overpartition pairs with conditions
on their Durfee dissection, as well as certain lattice paths.
When further specialized, the series become infinite products,
leading to numerous identities for partitions, overpartitions, and
overpartition pairs.
\end{abstract}

\maketitle

\section{Statement of Results}
In this paper we study two classes of basic hypergeometric series,
\begin{eqnarray} \label{rki}
R_{k,i}(a,b;x;q) &=&
\frac{(-axq,-bxq)_{\infty}}{(xq,abxq)_{\infty}} \sum_{n\geq 0}
\frac{(-ab)^nx^{kn}q^{kn^2 +(k-i+1)n - \binom{n}{2}}
(-1/a,-1/b)_n(xq)_n}{(q,-axq,-bxq)_n} \nonumber \\
&\times& \left( 1 - \frac{abx^iq^{(2n+1)i -
2n}(1+q^n/a)(1+q^n/b)}{(1+axq^{n+1})(1+bxq^{n+1})} \right)
\end{eqnarray}
and
\begin{eqnarray} \label{rkitilde}
\tilde{R}_{k,i}(a,b;x;q) &=&
\frac{(-axq,-bxq)_{\infty}}{(xq,abxq)_{\infty}} \sum_{n \geq 0}
\frac{(-ab)^nx^{(k-1)n}q^{kn^2 + (k -i)n -2\binom{n}{2}}(-1/a,-1/b)_n(x^2q^2;q^2)_n}{(q^2;q^2)_n(-axq,-bxq)_n} \nonumber \\
&\times& \left(1 - \frac{abx^iq^{(2n+1)i -
2n}(1+q^n/a)(1+q^n/b)}{(1+axq^{n+1})(1+bxq^{n+1})} \right).
\end{eqnarray}
Here we have employed the standard $q$-series notation
\cite{Ga-Ra1}
\begin{equation}
(a;q)_n = \prod_{j = 0}^{n-1} (1-aq^j)
\end{equation}
and
\begin{equation}
(a_1,a_2,...,a_k)_n = (a_1;q)_n(a_2;q)_n\cdots(a_k;q)_n.
\end{equation}

In the first part of the paper we interpret the coefficient of
$a^sb^tx^mq^n$ in \eqref{rki} and \eqref{rkitilde} in terms of
overpartition pairs. Recall that an overpartition is a partition
in which the first occurrence of a number may be overlined.  To
speak concisely about the relevant overpartition pairs, we shall
say that $j$ occurs \emph{unattached} in the overpartition pair
$(\lambda,\mu)$ if it only occurs non-overlined and only in $\mu$.
For example, in the overpartition pair
$((\overline{6},4,4,3),(6,\overline{4},4,\overline{2},2,1))$, only
$1$ occurs unattached.  We also define the \emph{valuation} of an
overpartition pair $(\lambda,\mu)$ at $j$ as
\begin{equation} \label{val}
v_j((\lambda,\mu)) = f_{j}(\lambda) + f_{\overline{j}}(\lambda) +
f_{\overline{j}}(\mu) + \chi(\text{$j$ occurs unattached in
$(\lambda,\mu)$}),
\end{equation}
where $\chi$ is the usual characteristic function and
$f_j(\lambda)$ counts the number of occurrences of $j$ in
$\lambda$.  We are now prepared to state our first two theorems.
Here and throughout the paper we assume that $k \geq 2$ and $1
\leq i \leq k$, unless otherwise noted.
\begin{theorem} \label{thm1}
If
$$
R_{k,i}(a,b;x;q) = \sum_{s,t,m,n \geq
0}r_{k,i}(s,t,m,n)a^sb^tx^mq^n,
$$
then $r_{k,i}(s,t,m,n)$ is equal to the number of overpartition
pairs $(\lambda,\mu)$ of $n$ with $m$ parts, $s$ of which are
overlined and in $\lambda$ or non-overlined and in $\mu$, $t$ of
which are in $\mu$, where $(i)$ $v_1((\lambda,\mu)) \leq i-1$, and
$(ii)$ for each $j \geq 1$, $f_j(\lambda) + v_{j+1}((\lambda,\mu))
\leq k-1$.
\end{theorem}

\begin{theorem} \label{thm2}
If
$$
\tilde{R}_{k,i}(a,b;x;q) = \sum_{s,t,m,n \geq
0}\tilde{r}_{k,i}(s,t,m,n)a^sb^tx^mq^n,
$$
then $\tilde{r}_{k,i}(s,t,m,n)$ is equal to the number of
overpartition pairs $(\lambda,\mu)$ counted by $r_{k,i}(s,t,m,n)$
such that if there is equality in condition $(ii)$ of Theorem
\ref{thm1} for some $j \geq 1$, then
\begin{equation} \label{three}
jf_j(\lambda) + (j+1)v_{j+1}((\lambda,\mu)) \equiv i-1 +
\mathcal{O}_j(\lambda) + \mathcal{O}_j(\mu) \pmod{2},
\end{equation}
where $\mathcal{O}_j(\cdot)$ denotes the number of overlined parts
less than or equal to $j$.
\end{theorem}

These theorems unify and generalize many important families of
partition identities, including Gordon's generalization of the
Rogers-Ramanujan identities \cite{Go1},  Bressoud's
Rogers-Ramanujan identities for even moduli \cite{Br1}, Gordon's
theorems for overpartitions \cite{Lo1},  Andrews' generalization
of the G\"ollnitz-Gordon identities \cite{An.1}, their
overpartition analogue \cite{Lo2}, as well as some more general
results of Corteel and the authors \cite{Co-Lo-Ma1,Co-Ma1}.  How
all of these results follow from Theorems \ref{thm1} and
\ref{thm2} will be explained in Section \ref{corsec}, and several
new families of identities will be presented.

In the second part of the paper we study three more classes of
combinatorial objects counted by $R_{k,i}(a,b;1;q)$ and
$\tilde{R}_{k,i}(a,b;1;q)$.  It will be necessary to defer the
definitions of these objects to later in the paper.

\begin{theorem} \label{thm3}
Let $B_{k,i}(s,t,n)$ denote the number of overpartition pairs
which are counted by $r_{k,i}(s,t,m,n)$ for some $m$. Let
$C_{k,i}(s,t,n)$ denote the number of overpartition pairs of $n$
whose Frobenius representations have $s$ non-overlined parts in
their bottom rows and $t$ non-overlined parts in their top rows,
and whose successive ranks are in the interval $[-i+2,2k-i-1]$.
Let $D_{k,i}(s,t,n)$ denote the number of $(k,i)$-admissible
overpartition pairs of $n$ whose Frobenius representations have
$s$ non-overlined parts in their bottom rows and $t$ non-overlined
parts in their top rows. Let $E_{k,i}(s,t,n)$ denote the number of
generalized Bressoud-Burge lattice paths of major index $n$
satisfying the odd $(k,i)$-conditions, where the number of peaks
marked by $a$ (resp.\ marked by $b$) is $s$ (resp.\ $t$). Then
$$
B_{k,i}(s,t,n) = C_{k,i}(s,t,n) = D_{k,i}(s,t,n) = E_{k,i}(s,t,n).
$$
\end{theorem}

\begin{theorem} \label{thm4}
Let $\tilde{B}_{k,i}(s,t,n)$ denote the number of overpartition
pairs which are counted by $\tilde{r}_{k,i}(s,t,m,n)$ for some
$m$. Let $\tilde{C}_{k,i}(s,t,n)$ denote the number of
overpartition pairs of $n$ whose Frobenius representations have
$s$ non-overlined parts in their bottom rows and $t$ non-overlined
parts in their top rows, and whose successive ranks are in the
interval $[-i+2,2k-i-2]$. Let $\tilde{D}_{k,i}(s,t,n)$ denote the
number of self-$(k,i)$-conjugate overpartition pairs of $n$ whose
Frobenius representations have $s$ non-overlined parts in their
bottom rows and $t$ non-overlined parts in their top rows.  Let
$\tilde{E}_{k,i}(s,t,n)$ denote the number of generalized
Bressoud-Burge lattice paths counted by $E_{k,i}(s,t,n)$ which
also satisfy the even $(k,i)$-conditions.  Then
$$
\tilde{B}_{k,i}(s,t,n) = \tilde{C}_{k,i}(s,t,n) =
\tilde{D}_{k,i}(s,t,n) = \tilde{E}_{k,i}(s,t,n).
$$
\end{theorem}

Theorems \ref{thm3} and \ref{thm4} extend overpartition-theoretic
work of Corteel and the authors \cite{Co-Lo-Ma1,Co-Ma1}, which had
in turn generalized partition-theoretic work of Andrews, Bressoud,
and Burge \cite{An2.5,An5,Br2.5,Br3,Bu1,Bu2}.

The paper is organized as follows: Theorem \ref{thm1} is proven in
Section $2$ using Andrews' $q$-difference equations for some
families of basic hypergeometric series \cite{An1}. Theorem
\ref{thm2} is proven in Section $3$ in the same way, except that
we will have to develop the required $q$-difference equations from
scratch.  In Section $4$, we present some of the many
combinatorial identities which follow from Theorems \ref{thm1} and
\ref{thm2}.  In Sections $5 - 7$ we define the combinatorial
structures occurring in Theorems \ref{thm3} and \ref{thm4} and
prove these theorems.

\section{The $R_{k,i}(a,b;x;q)$}

It was Andrews who first observed the combinatorial significance
of series like \eqref{rki} and \eqref{rkitilde}.  Selberg
\cite{Se1} had essentially proven $q$-difference equations for
$R_{k,i}(0,0;x;q)$, and Andrews \cite{An.01} showed how this could
be used to prove Gordon's generalization of the Rogers-Ramanujan
identities \cite{Go1}.  He then proceeded to develop a massive
generalization of the $R_{k,i}(0,0;x;q)$ \cite{An1}, and the
combinatorics of these series has turned out to be one of the
major areas of research in the theory of partitions over the last
$40$ years (e.g.
\cite{An.05,An.06,An.1,An2,An2.5,An3,An5,An-Sa1,Br1,Br2,Br3,Bu1,Bu2,An-Br1,Con1,Co-Ma1,Co-Lo-Ma1,Lo1,Lo2,Lo3,Lo4}).
The series also have direct applications to $q$-series identities
(e.g. \cite{An1,An4}) and $q$-continued fractions (e.g.
\cite{An1,An-Be1}).

In terms of Andrews' series, called
$J_{\lambda,k,i}(a_1,a_2,\dots,a_{\lambda};x;q)$ \cite{An1}, we
have
\begin{equation}
R_{k,i}(a,b;x;q) =
\frac{1}{(abxq)_{\infty}}J_{2,k,i}(-1/a,-1/b;x;q).
\end{equation}
Employing the $q$-difference equations for the $J_{2,k,i}$ (and
related functions) \cite[Eq. (2.1)-(2.4)]{An1}, we may deduce that
the $R_{k,i}(a,b;x;q)$ satisfy the following:
\begin{lemma} \label{rkirecurlem}
\begin{equation} \label{rkirecur1}
R_{k,1}(a,b;x;q) = R_{k,k}(a,b;xq;q),
\end{equation}
\begin{eqnarray} \label{rkirecur2}
R_{k,2}(a,b;x;q) - R_{k,1}(a,b;x;q) &=&
\frac{xq}{(1-abxq)}R_{k,k-1}(a,b;xq;q)  \\ &+&
\frac{axq}{(1-abxq)}R_{k,k}(a,b;xq;q) \nonumber \\ &+&
\frac{bxq}{(1-abxq)}R_{k,k}(a,b;xq;q) \nonumber \\ &+&
\frac{abxq}{(1-abxq)}R_{k,k}(a,b;xq;q) \nonumber,
\end{eqnarray}
and for $3 \leq i \leq k$,
\begin{eqnarray} \label{rkirecur3}
R_{k,i}(a,b;x;q) - R_{k,i-1}(a,b;x;q) &=&
\frac{(xq)^{i-1}}{(1-abxq)} R_{k,k-i+1}(a,b;xq;q) \\ &+&
\frac{a(xq)^{i-1}}{(1-abxq)}R_{k,k-i+2}(a,b;xq;q) \nonumber \\ &+&
\frac{b(xq)^{i-1}}{(1-abxq)} R_{k,k-i+2}(a,b;xq;q) \nonumber \\
&+& \frac{ab(xq)^{i-1}}{(1-abxq)}R_{k,k-i+3}(a,b;xq;q)) \nonumber.
\end{eqnarray}
\end{lemma}
Using these, we may deduce Theorem \ref{thm1}.
\newline \emph{Proof of Theorem \ref{thm1}.}
First, observe that the $q$-difference equations in Lemma
\ref{rkirecurlem} together with the fact that $R_{k,i}(a,b;0;q) =
1$ uniquely define the functions $R_{k,i}(a,b;x;q)$.   Now, let
$$
\widehat{R}_{k,i}(a,b;x;q) = \sum_{s,t,m,n \geq
0}r_{k,i}(s,t,m,n)a^sb^tx^mq^n.
$$
We wish to show that $\widehat{R}_{k,i}(a,b;x;q) =
R_{k,i}(a,b;x;q)$.  We shall accomplish this by showing that the
functions $\widehat{R}_{k,i}(a,b;x;q)$ satisfy the same
$q$-difference equations as the $R_{k,i}(a,b;x;q)$ in Lemma
\ref{rkirecurlem}. The fact that $\widehat{R}_{k,i}(a,b;0;q) = 1$
is obvious, since there are no overpartition pairs without any
parts except for the empty one.

Observe that subtracting one from each part of an overpartition
pair (and deleting the resulting zeros) that satisfies condition
$(ii)$ in Theorem \ref{thm1} gives another overpartition pair that
satisfies this condition. Similarly, adding one to each part of an
overpartition pair that satisfies the condition gives another
overpartition pair that satisfies the condition.

We begin with \eqref{rkirecur1}.  An overpartition pair
$(\lambda,\mu)$ counted by $r_{k,1}(s,t,m,n)$ has no ones
whatsoever and hence has $v_2((\lambda,\mu)) \leq k-1$.  By
subtracting one from each part we see that
$\widehat{R}_{k,1}(a,b;x;q) = \widehat{R}_{k,k}(a,b;xq;q)$.

For \eqref{rkirecur2}, we observe that the function
$$
\widehat{R}_{k,2}(a,b;x;q) - \widehat{R}_{k,1}(a,b;x;q)
$$
is the generating function for those overpartition pairs
$(\lambda,\mu)$ counted by $r_{k,2}(s,t,m,n)$ having
$v_1((\lambda,\mu)) = 1$.   We break these pairs into four
disjoint classes:  those having $1$ as a part of $\lambda$, those
having $\overline{1}$ as a part of $\lambda$, those having
$\overline{1}$ as a part of $\mu$, and those in which $1$ occurs
unattached.   In the first of these four cases,
$v_2((\lambda,\mu)) \leq k-2$.  So, removing the $1$ from
$\lambda$ along with any ones that may occur in $\mu$, and then
subtracting one from all of the remaining parts, we see that these
overpartition pairs are generated by
$$
\frac{xq}{(1-abxq)}\widehat{R}_{k,k-1}(a,b;xq;q).
$$
In the second of these cases, where $\overline{1}$ occurs in
$\lambda$,  we have $v_2((\lambda,\mu)) \leq k-1$.  So,  removing
the $\overline{1}$ from $\lambda$ along with any ones that may
occur in $\mu$, and then subtracting one from all of the remaining
parts, we see that these overpartition pairs are generated by
$$
\frac{axq}{(1-abxq)}\widehat{R}_{k,k}(a,b;xq;q).
$$
In exactly the same way we see that those pairs containing an
$\overline{1}$ in $\mu$ are generated by
$$
\frac{bxq}{(1-abxq)}\widehat{R}_{k,k}(a,b;xq;q).
$$
For the final case, where $1$ occurs unattached in the
overpartition pair $(\lambda,\mu)$, again we have
$v_2((\lambda,\mu) \leq k-1$.  So, removing all of the ones from
$\mu$ and subtracting one from all of the remaining parts, we see
that these overpartition pairs are generated by
$$
\frac{abxq}{(1-abxq)}\widehat{R}_{k,k}(a,b;xq;q).
$$
Putting everything together gives \eqref{rkirecur2}.

We now turn to \eqref{rkirecur3}.   As above, the function
$$
\widehat{R}_{k,i}(a,b;x;q) - \widehat{R}_{k,i-1}(a,b;x;q)
$$
is the generating function for those overpartition pairs which are counted
by $r_{k,i}(s,t,m,n)$ and which have $v_1((\lambda,\mu)) = i-1$.   And, as
before, we consider four cases:  $f_1(\lambda) = i-1$,
$f_1(\lambda) = i-2$ and $f_{\overline{1}}(\lambda) = 1$,
$f_1(\lambda) = i-2$ and $f_{\overline{1}}(\mu) = 1$, and
$f_1(\lambda) = i-3$, $f_{\overline{1}}(\lambda) = 1$, and
$f_{\overline{1}}(\mu) = 1$.   Notice that since $i-1 \geq 2$ we
cannot have an unattached occurrence of $1$ in $(\lambda,\mu)$.
Now, in the first of these cases, $v_2((\lambda,\mu)) \leq k-i$.
So, removing the $i-1$ ones from $\lambda$ as well as any
non-overlined ones from $\mu$, and then subtracting one from each
remaining part, we see that these overpartition pairs are
generated by
$$
\frac{(xq)^{i-1}}{(1-abxq)}\widehat{R}_{k,k-i+1}(a,b;xq;q).
$$
For the second case, where $f_1(\lambda) = i-2$ and
$f_{\overline{1}}(\lambda) = 1$, we have $v_2((\lambda,\mu)) \leq
k-i+1$.  So, removing the $i-2$ ones and the $\overline{1}$ from
$\lambda$, as well as any non-overlined ones from $\mu$, and then
subtracting one from each remaining part, we see that these
overpartition pairs are generated by
$$
\frac{axq(xq)^{i-2}}{(1-abxq)}\widehat{R}_{k,k-i+2}(a,b;xq;q).
$$
Similarly, those overpartition pairs having $f_1(\lambda) = i-2$
and $f_{\overline{1}}(\mu)$ are generated by
$$
\frac{bxq(xq)^{i-2}}{(1-abxq)}\widehat{R}_{k,k-i+2}(a,b;xq;q).
$$
Finally, if $f_1(\lambda) = i-3$, $f_{\overline{1}}(\lambda) = 1$,
and $f_{\overline{1}}(\mu) = 1$, then $v_2((\lambda,\mu) \leq
k-i+2$.  So, removing the $i-3$ ones and the $\overline{1}$ from
$\lambda$, the $\overline{1}$ and any non-overlined ones from
$\mu$, and then subtracting one from each remaining part, we see
that these overpartition pairs are generated by
$$
\frac{(axq)(bxq)(xq)^{i-3}}{(1-abxq)}\widehat{R}_{k,k-i+3}(a,b;xq;q).
$$
Putting everything together gives \eqref{rkirecur3} for the
$\widehat{R}_{k,i}(a,b;x;q)$ and we may now conclude that
$R_{k,i}(a,b;x;q) = \widehat{R}_{k,i}(a,b;x;q)$, establishing
Theorem \ref{thm1}. \qed

\section{The $\tilde{R}_{k,i}(a,b;x;q)$}

Unlike the case for the $R_{k,i}(a,b;x;q)$, we will need to
develop from scratch the theory of recurrences for the
$\tilde{R}_{k,i}(a,b;x;q)$.  In this endeavor we closely follow
Andrews \cite{An1}.  For $k \geq 1$ and $i \in \mathbb{Z}$, define
\small$$ 
\tilde{H}_{2,k,i}(a,b;x;q) =
\frac{(-axq,-bxq)_{\infty}}{(xq)_{\infty}} \sum_{n \geq 0}
\frac{(-ab)^nx^{(k-1)n}q^{kn^2 +n -in
-2\binom{n}{2}}(-1/a,-1/b)_n(x^2;q^2)_n(1-x^iq^{2ni})}{(q^2;q^2)_n(-axq,-bxq)_n(1-x)}
$$
\normalsize
and
$$
\tilde{J}_{2,k,i}(a,b;x;q) =
(abxq)_{\infty}\tilde{R}_{k,i}(a,b;x;q).
$$
When $a=0$ or $b=0$, these functions simplify to the
$\tilde{H}_{k,i}$ and $\tilde{J}_{k,i}$ studied in
\cite{Co-Lo-Ma1}.  We shall establish the following facts:

\begin{proposition} \label{htilderecur}
We have
\begin{equation} \label{htilde1}
\tilde{H}_{2,k,0}(a,b;x;q) = 0,
\end{equation}
\begin{equation} \label{htilde2}
\tilde{H}_{2,k,-i}(a,b;x;q) = -x^{-i}\tilde{H}_{2,k,i}(a,b;x;q),
\end{equation}
\begin{equation} \label{htilde3}
\tilde{H}_{2,k,i}(a,b;x;q) - \tilde{H}_{2,k,i-2}(a,b;x;q) =
x^{i-2}(1+x)\tilde{J}_{2,k,k-i+1}(a,b;x;q),
\end{equation}
and
\begin{equation} \label{htilde4}
\tilde{J}_{2,k,i}(a,b;x;q) = \tilde{H}_{2,k,i}(a,b;xq;q) +
(axq+bxq)\tilde{H}_{2,k,i-1}(a,b;xq;q) +
abx^2q^2\tilde{H}_{2,k,i-2}(a,b;xq;q).
\end{equation}
\end{proposition}

\begin{proof}
Equations \eqref{htilde1} and \eqref{htilde2} are straightforward.
For the other two, we introduce a little notation to simplify the
calculations.  We write
\begin{equation}
\tilde{C}_{2,k,i}(a,b;x;q) =
\frac{(xq)_{\infty}}{(-axq,-bxq)_{\infty}}
\tilde{H}_{2,k,i}(a,b;x;q)
\end{equation}
and
\begin{equation}
\tilde{D}_{2,k,i}(a,b;x;q) =
\frac{(xq)_{\infty}}{(-axq,-bxq)_{\infty}}
\tilde{J}_{2,k,i}(a,b;x;q).
\end{equation}
We also write
\begin{equation}
\tilde{M}_n(a,b;x;q) = \tilde{M}_n(x) = x^{(k-1)n}q^{(k-1)n^2 +
2n}(-ab)^n.
\end{equation}
These three definitions are in analogy with Andrews' definitions
in \cite{An1}.   We note that
\begin{equation} \label{mtilde1}
\tilde{M}_{n+1}(x) = -abx^{k-1}q^{2n(k-1) + k +1} \tilde{M}_n(x)
\end{equation}
and
\begin{equation} \label{mtilde2}
\tilde{M}_n(xq) = q^{(k-1)n}\tilde{M}_n(x).
\end{equation}

We are now prepared to deal with \eqref{htilde3}.  We have
\begin{eqnarray*}
\tilde{C}_{2,k,i}(a,b;x;q) &-& \tilde{C}_{2,k,i-2}(a,b;x;q) \\ &=&
\sum_{n \geq 0} \frac{\tilde{M}_n(x)(x^2;q^2)_n(-1/a,-1/b)_n}
{(1-x)(q^2;q^2)_n(-axq,-bxq)_n} \\
&\times& \left(q^{-in}(1-q^{2n}) + x^{i-2}q^{n(i-2)}(1-x^2q^{2n})
\right) \\
&=& (1+x) \sum_{n \geq 1}
\frac{\tilde{M}_n(x)q^{-in}(x^2q^2;q^2)_{n-1}(-1/a,-1/b)_n}{(q^2;q^2)_{n-1}(-axq,-bxq)_n}
\\
&+& x^{i-2}(1+x) \sum_{n \geq 0}
\frac{\tilde{M}_{n}(x)q^{n(i-2)}(x^2q^2;q^2)_n(-1/a,-1/b)_n}{(q^2;q^2)_n(-axq,-bxq)_n}
\\
&=& (1+x)\sum_{n \geq 0}
\frac{\tilde{M}_{n+1}(x)q^{-i(n+1)}(x^2q^2;q^2)_{n}(-1/a,-1/b)_{n+1}}{(q^2;q^2)_{n}(-axq,-bxq)_{n+1}}
\\
&+& x^{i-2}(1+x) \sum_{n \geq 0}
\frac{\tilde{M}_{n}(x)q^{n(i-2)}(x^2q^2;q^2)_n(-1/a,-1/b)_n}{(q^2;q^2)_n(-axq,-bxq)_n}
\\
&=& -abx^{k-1}q^{k-i+1}(1+x)\sum_{n \geq
0}\frac{\tilde{M}_n(x)q^{n(2k-i-2)}(x^2q^2;q^2)_n(-1/a,-1/b)_{n+1}}{(q^2;q^2)_n(-axq,-bxq)_{n+1}}
\\
&+& x^{i-2}(1+x) \sum_{n \geq 0}
\frac{\tilde{M}_{n}(x)q^{n(i-2)}(x^2q^2;q^2)_n(-1/a,-1/b)_n}{(q^2;q^2)_n(-axq,-bxq)_n}
\\
&=& -abx^{k-1}q^{k-i+1}(1+x) \sum_{n \geq 0}
\frac{\tilde{M}_n(xq)q^{n(k-i-1)}
(x^2q^2;q^2)_n(-1/a,-1/b)_{n+1}}{(q^2;q^2)_n(-axq,-bxq)_{n+1}} \\
&+& x^{i-2}(1+x)\sum_{n \geq 0}
\frac{\tilde{M}_n(xq)q^{n(-k+i-1)}(x^2q^2;q^2)_n(-1/a,-1/b)_n}{(q^2;q^2)_n(-axq,-bxq)_n}
\\
&=& x^{i-2}(1+x)\sum_{n \geq 0}
\frac{\tilde{M}_n(xq)q^{-n(k-i+1)}(x^2q^2;q^2)_n(-1/a,-1/b)_n}{(q^2;q^2)_n(-axq,-bxq)_n}
\\
&-& abx^{i-2}(1+x)(xq)^{k-i+1}\\
&&\quad\quad \times \sum_{n \geq 0}
\frac{\tilde{M}_n(xq)q^{-n(k-i+1) + 2n(k-i+1) -
2n}(x^2q^2;q^2)_n(-1/a,-1/b)_{n+1}}{(q^2;q^2)_n(-axq,-bxq)_{n+1}}
\\
&=& x^{i-2}(1+x)\sum_{n \geq 0} \frac{\tilde{M}_n(xq)q^{-n(k-i+1)}
(x^2q^2;q^2)_n(-1/a,-1/b)_n}{(q^2;q^2)_n(-axq,-bxq)_n} \\
&\times& \left( 1 - \frac{abx^{k-i+1}q^{(2n+1)(k-i+1) -
2n}(1+q^n/a)(1+q^n/b)}{(1-axq^{n+1})(1-bxq^{n+1})} \right) \\
&=& x^{i-2}(1+x)\tilde{D}_{2,k,k-i+1}(a,b;x;q).
\end{eqnarray*}
Multiplying the extremes of the above string of equations by
$(-axq,-bxq)_{\infty}/(xq)_{\infty}$ yields \eqref{htilde3}.

We now turn to \eqref{htilde4}.   By making a common denominator
in the expression in parentheses in the definition of the
$\tilde{J}_{2,k,i}(a,b;x;q)$, we have
\small 
\begin{eqnarray*}
\tilde{D}_{2,k,i}(a,b;x;q) &=& \sum_{n \geq 0}
\frac{\tilde{M}_n(xq)q^{-in}(x^2q^2;q^2)_n(-1/a,-1/b)_n}{(q^2;q^2)_n(-axq,-bxq)_{n+1}}
\\
&\times& \left( 1 + (a+b)xq^{n+1} + abx^2q^{2n+2} -
x^iq^{(2n+1)i}(1 + (a+b)q^{-n} + abq^{-2n}) \right) \\
&=& \sum_{n \geq 0}
\frac{\tilde{M}_n(xq)q^{-in}(x^2q^2;q^2)_n(-1/a,-1/b)_n}{(q^2;q^2)_n(-axq,-bxq)_{n+1}}
\\
&\times& \left[(1-x^iq^{(2n+1)i}) +
(a+b)xq^{n+1}(1-x^{i-1}q^{(2n+1)(i-1)}) \right. \\
&+& \left. abx^2q^{2n+2}(1 - x^{i-2}q^{(2n+1)(i-2)})\right]
\\
&=& \frac{(1-xq)}{(1+axq)(1+bxq)} \left[ \sum_{n \geq 0}
\frac{\tilde{M}_n(xq)q^{-in}((xq)^2;q^2)_n(-1/a,-1/b)_n(1-(xq)^iq^{2ni})}
{(q^2;q^2)_n(-a(xq)q,-b(xq)q)_n(1-xq)} \right.\\
&+& (a+b)xq \sum_{n \geq 0}
\frac{\tilde{M}_n(xq)q^{-(i-1)n}((xq)^2;q^2)_n(-1/a,-1/b)_n(1-(xq)^{i-1}q^{2n(i-1)})}
{(q^2;q^2)_n(-a(xq)q,-b(xq)q)_n(1-xq)} \\
&+& \left. abx^2q^2 \sum_{n \geq 0}
\frac{\tilde{M}_n(xq)q^{-(i-2)n}((xq)^2;q^2)_n(-1/a,-1/b)_n(1-(xq)^{i-2}q^{2n(i-2)})}
{(q^2;q^2)_n(-a(xq)q,-b(xq)q)_n(1-xq)} \right] \\
&=& \frac{(1-xq)}{(1+axq)(1+bxq)} \left[ \tilde{C}_{2,k,i}(a,b;xq;q) \right.\\
&+& \left. (a+b)xq\tilde{C}_{2,k,i-1}(a,b;xq;q) +
abx^2q^2\tilde{C}_{2,k,i-2}(a,b;xq;q) \right].
\end{eqnarray*}
\normalsize Multiplying both sides of this string of equations by
$(-axq,-bxq)_{\infty}/(xq)_{\infty}$ finishes the proof of
\eqref{htilde4}.  And this then completes the proof of Proposition
\ref{htilderecur}.
\end{proof}

We now have the analogue of Lemma \ref{rkirecurlem} for the
$\tilde{R}_{k,i}(a,b;x;q)$ using Proposition \ref{htilderecur}.
\begin{lemma} \label{rkitilderecurlem}
\begin{equation} \label{rkitilderecur1}
\tilde{R}_{k,1}(a,b;x;q) = \tilde{R}_{k,k}(a,b;xq;q),
\end{equation}
\begin{eqnarray} \label{rkitilderecur2}
\tilde{R}_{k,2}(a,b;x;q) &=&
\frac{1}{(1-abxq)}\tilde{R}_{k,k-1}(a,b;xq;q) \\ &+&
\frac{xq}{(1-abxq)}\tilde{R}_{k,k-1}(a,b;xq;q) \nonumber \\ &+&
\frac{axq}{(1-abxq)} \tilde{R}_{k,k}(a,b;xq;q) \nonumber \\ &+&
\frac{bxq}{(1-abxq)} \tilde{R}_{k,k}(a,b;xq;q) \nonumber,
\end{eqnarray}
and, for $3 \leq i \leq k$,
\begin{eqnarray} \label{rkitilderecur3}
\tilde{R}_{k,i}(a,b;x;q) - \tilde{R}_{k,i-2}(a,b;x;q) &=&
\frac{(xq)^{i-2}(1+xq)}{(1-abxq)} \tilde{R}_{k,k-i+1}(a,b;xq;q) \\
&+& \frac{a(xq)^{i-2}(1+xq)}{(1-abxq)}
\tilde{R}_{k,k-i+2}(a,b;xq;q) \nonumber \\ &+&
\frac{b(xq)^{i-2}(1+xq)}{(1-abxq)} \tilde{R}_{k,k-i+2}(a,b;xq;q)
\nonumber \\ &+&  \frac{ab(xq)^{i-2}(1+xq)}{(1-abxq)}
\tilde{R}_{k,k-i+3}(a,b;xq;q) \nonumber.
\end{eqnarray}
\end{lemma}

\emph{Proof of Theorem \ref{thm2}.} We begin by observing that for
an overpartition pair $(\lambda,\mu)$ counted by
$\tilde{r}_{k,i}(s,t,m,n)$, if $(\lambda,\mu) - \vec{1}$ satisfies
condition $(ii)$ of Theorem \ref{thm1} at $j$, then so does
$(\lambda,\mu)$ at $j+1$, where $- \vec{1}$ is shorthand for
subtracting $1$ from each part and then deleting any zeros. Hence,
we have
\begin{eqnarray} \label{special}
jf_j(\lambda - \vec{1}) + (j+1)v_{j+1}((\lambda,\mu) - \vec{1})
& = & (j+1)f_{j+1}(\lambda) + (j+2)v_{j+2}((\lambda,\mu))  \nonumber \\
&-&  \left(f_{j+1}(\lambda) + v_{j+2}((\lambda,\mu))\right)  \nonumber \\
&\equiv& i-1 + \mathcal{O}_{j+1}(\lambda) + \mathcal{O}_{j+1}(\mu)
- (k-1) \pmod{2} \nonumber \\
&\equiv& k-i + \mathcal{O}_j(\lambda-\vec{1}) + \mathcal{O}_j(\mu
- \vec{1})  \\
&+& \begin{cases} 0, & \text{if $\overline{1} \not \in
\lambda$ and $\overline{1} \not \in \mu$}, \\
0, & \text{if $\overline{1} \in
\lambda$ and $\overline{1} \in \mu$}, \\
1, & \text{otherwise}
\end{cases}
\pmod{2}. \nonumber
\end{eqnarray}

We now proceed as in the proof of Theorem \ref{thm1}.  Equation
\eqref{special} will be used throughout the proof to ensure
condition \eqref{three} in the overpartition pairs under
consideration.  This may not always be mentioned explicitly.  We
observe that since $\tilde{R}_{k,i}(a,b;0;q) = 1$, the
$q$-difference equations in Lemma \ref{rkitilderecurlem} uniquely
define the functions $\tilde{R}_{k,i}(a,b;x;q)$.  Now let
$$
\tilde{S}_{k,i}(a,b;x;q) = \sum_{s,t,m,n \geq 0}
\tilde{r}_{k,i}(s,t,m,n)a^sb^tx^mq^n.
$$
Again, $\tilde{S}_{k,i}(a,b;0;q) = 1$ because there is only one
overpartition pair without parts - the empty one.  So, to prove
Theorem \ref{thm2} we need to show that the
$\tilde{S}_{k,i}(a,b;x;q)$ satisfy the same $q$-difference
equations as the $\tilde{R}_{k,i}(a,b;x;q)$ in Lemma
\ref{rkitilderecurlem}.

We begin with \eqref{rkitilderecur1}.  An overpartition pair
counted by $\tilde{r}_{k,1}(s,t,m,n)$ has no ones and has
$v_2((\lambda,\mu)) \leq k-1$.  Subtracting one from each part and
appealing to \eqref{special}, we see that these overpartition
pairs are generated by $\tilde{S}_{k,k}(a,b;xq;q)$.

For \eqref{rkitilderecur2}, an overpartition pair $(\lambda,\mu)$
counted by $\tilde{r}_{k,2}(s,t,m,n)$ has either no ones or
$v_1((\lambda,\mu))=1$.   If there are no ones, then
$v_2((\lambda,\mu)) \leq k-2$.  Notice that in this case we cannot
have $v_2((\lambda,\mu)) = k-1$, for then we would have
$1f_1((\lambda,\mu)) + 2v_2((\lambda,\mu)) \equiv 0 \pmod{2}$,
which violates the condition \eqref{three} defining the
$\tilde{r}_{k,2}(s,t,m,n)$.  Hence we have $v_2((\lambda,\mu) \leq
k-2$.  So, subtracting one from each part of $(\lambda,\mu)$ and
appealing to \eqref{special}, we see that these pairs are
generated by
\begin{equation} \label{s1}
\tilde{S}_{k,k-1}(a,b;xq;q)
\end{equation}

Now, if $v_1((\lambda,\mu)) = 1$, this may be for one of four
reasons: $1$ occurs in $\lambda$, $\overline{1}$ occurs in
$\lambda$, $\overline{1}$ occurs in $\mu$, or $1$ occurs
unattached.  In the first case, we have $v_2((\lambda,\mu)) \leq
k-2$.  Subtracting one from each part, removing the $1$ from
$\lambda$ as well as any non-overlined ones from $\mu$, these
pairs are seen to be generated by
\begin{equation} \label{s2}
\frac{xq}{(1-abxq)}\tilde{S}_{k,k-1}(a,b;xq;q).
\end{equation}
In the second case, when $\overline{1}$ appears in $\lambda$, then
we have $v_2((\lambda,\mu)) \leq k-1$.   Notice that if
$v_2((\lambda,\mu)) = k-1$, then we have $1f_1(\lambda) +
2v_2((\lambda,\mu)) \equiv 0 \pmod{2}$, which is congruent to $2-1
+ \mathcal{O}_1(\lambda) + \mathcal{O}_1(\mu)$ modulo $2$.
Removing the $\overline{1}$ from $\mu$, removing any non-overlined
ones from $\mu$, and then subtracting one from each remaining
part, we find (keeping in mind \eqref{special}) that these pairs
are generated by
\begin{equation} \label{s3}
\frac{axq}{(1-abxq)}\tilde{S}_{k,k}(a,b;xq;q).
\end{equation}
The third case, when $\overline{1}$ occurs in $\mu$, is analogous
to the second case and these overpartition pairs are generated by
\begin{equation} \label{s4}
\frac{bxq}{(1-abxq)}\tilde{S}_{k,k}(a,b;xq;q).
\end{equation}
Finally, we consider the case when $1$ occurs unattached in the
overpartition pair $((\lambda,\mu))$.  If $v_2((\lambda,\mu)) =
k-1$, then condition \eqref{three} in the definition of the
$r_{k,2}(s,t,m,n)$ would be violated, so we have
$v_2((\lambda,\mu)) \leq k-2$.  Removing all of the unattached
ones and subtracting one from each remaining part, we see that
these pairs are generated by
\begin{equation} \label{s5}
\frac{abxq}{(1-abxq)}\tilde{S}_{k,k-1}(a,b;xq;q).
\end{equation}
Adding \eqref{s1} - \eqref{s5} together now shows that the
recurrence \eqref{rkitilderecur2} is true for the
$\tilde{S}_{k,i}(a,b;x;q)$.

We now turn to \eqref{rkitilderecur3}, proceeding much like
before.  The function $\tilde{S}_{k,i}(a,b;x;q) -
\tilde{S}_{k,i-2}(a,b;x;q)$ is the generating function for those
overpartition pairs $(\lambda,\mu)$ which are counted by
$\tilde{r}_{k,i}(s,t,m,n)$ and which have either
$v_1((\lambda,\mu)) = i-1$ or $v_1((\lambda,\mu)) = i-2$.  We
shall consider eight cases, the last of which has two subcases
depending on whether $i
> 3$.

In the first case, suppose that $v_1((\lambda,\mu)) = i-1$ and
$f_1(\lambda) = i-1$.  Then, $v_2((\lambda,\mu))$ can be as much
as $k-i$.  These pairs are generated by
\begin{equation} \label{t1}
\frac{(xq)^{i-1}}{(1-abxq)}\tilde{S}_{k,k-i+1}(a,b;xq;q).
\end{equation}
In the second case, suppose that $v_1((\lambda,\mu)) = i-2$ and
$f_1(\lambda) = i-2$.  Then, $v_2((\lambda,\mu))$ can be as much
as $k-i$, for if it were $k-i+1$ this would violate the condition
\eqref{three}.  These pairs are generated by
\begin{equation} \label{t2}
\frac{(xq)^{i-2}}{(1-abxq)}\tilde{S}_{k,k-i+1}(a,b;xq;q).
\end{equation}
In the third case, suppose that $v_1((\lambda,\mu)) = i-1$,
$f_{\overline{1}}(\lambda) = 1$, and $f_1(\lambda) = i-2$.  Then,
$v_2((\lambda,\mu))$ can be as much as $k-i+1$. These pairs are
generated by
\begin{equation} \label{t3}
\frac{a(xq)^{i-1}}{(1-abxq)}\tilde{S}_{k,k-i+2}(a,b;xq;q).
\end{equation}
In the fourth case, suppose that $v_1((\lambda,\mu)) = i-2$,
$f_{\overline{1}}(\lambda) = 1$, and $f_1(\lambda) = i-3$.  Then,
$v_2((\lambda,\mu))$ can be as much as $k-i+1$, for if it were
$k-i+2$ this would violate the condition \eqref{three}.  These
pairs are generated by
\begin{equation} \label{t4}
\frac{a(xq)^{i-2}}{(1-abxq)}\tilde{S}_{k,k-i+2}(a,b;xq;q).
\end{equation}
The fifth and sixth cases are analogous to the third and fourth,
respectively, where $f_{\overline{1}}(\lambda) = 1$ is replaced by
$f_{\overline{1}}(\mu) = 1$.   These pairs are generated by
\begin{equation} \label{t5}
\frac{b(xq)^{i-1} +
b(xq)^{i-2}}{(1-abxq)}\tilde{S}_{k,k-i+2}(a,b;xq;q).
\end{equation}
In the seventh case, suppose that $v_1((\lambda,\mu)) = i-1$,
$f_{\overline{1}}(\lambda) = f_{\overline{1}}(\mu) = 1$, and
$f_{1}(\lambda) = i-3$.  Then $v_2((\lambda,\mu)$ could be as much
as $k-i+2$.  These pairs are generated by
\begin{equation} \label{t6}
\frac{ab(xq)^{i-1}}{(1-abxq)}\tilde{S}_{k,k-i+3}(a,b;xq;q).
\end{equation}
For the eighth case, suppose that $v_1((\lambda,\mu)) = i-2$,
$f_{\overline{1}}(\lambda) = f_{\overline{1}}(\mu) = 1$, and
$f_{1}(\lambda) = i-4$.  This requires the assumption that $i \geq
4$.  Here $v_2((\lambda,\mu)) \leq k-i+2$, and these pairs are
generated (for $i \geq 4$) by
\begin{equation} \label{t7}
\frac{ab(xq)^{i-2}}{(1-abxq)}\tilde{S}_{k,k-i+3}(a,b;xq;q).
\end{equation}
Now, if $i = 3$ we cannot have $v_1((\lambda,\mu)) = i-2 = 1$
while at the same time having $\overline{1}$ occurring in both
$\lambda$ and $\mu$. What we can have, however, is $1$ occurring
unattached. Then $v_2((\lambda,\mu)) \leq k-1$ $(=k-i+2)$, and so
these pairs are generated by \eqref{t7} when $i=3$.

Adding together \eqref{t1} - \eqref{t7} establishes
\eqref{rkitilderecur3} for the $\tilde{S}_{k,i}(a,b;x;q)$ and we
may now conclude that $\tilde{S}_{k,i}(a,b;x;q) =
\tilde{R}_{k,i}(a,b;x;q)$, finishing the proof of Theorem
\ref{thm2}. \qed

\section{Corollaries} \label{corsec}
Using the fact that
$$
(a)_{-n} = \frac{(-1)^nq^{\binom{n+1}{2}}}{a^n(q/a)_n},
$$
the following representations for the $R_{k,i}(a,b;1;q)$ and
$\tilde{R}_{k,i}(a,b;1;q)$ can be deduced from \eqref{rki} and
\eqref{rkitilde}:
\begin{equation} \label{rkibilat}
R_{k,i}(a,b;1;q) = \frac{(-aq,-bq)_{\infty}}{(q,abq)_{\infty}}
\sum_{n \in \mathbb{Z}} \frac{q^{kn^2 + (k-i+1)n -
\binom{n}{2}}(-1/a,-1/b)_n(-ab)^n}{(-aq,-bq)_n}
\end{equation}
and
\begin{equation} \label{rkitildebilat}
\tilde{R}_{k,i}(a,b;1;q) =
\frac{(-aq,-bq)_{\infty}}{(q,abq)_{\infty}} \sum_{n \in
\mathbb{Z}} \frac{q^{kn^2+ (k-i)n -
2\binom{n}{2}}(-1/a,-1/b)_n(-ab)^n}{(-aq,-bq)_n}.
\end{equation}
Applying Jacobi's triple product identity,
\begin{equation} \label{jtp}
\sum_{n \in \mathbb{Z}} z^nq^{n^2} = (-zq,-q/z,q^2;q^2)_{\infty},
\end{equation}
one finds that many specializations of \eqref{rkibilat} and
\eqref{rkitildebilat} are infinite products with nice
combinatorial interpretations. Using Theorem \ref{thm1}, this
leads to numerous identities for partitions, overpartitions, and
overpartition pairs. For example, when $(a,b,q) = (1,1/q,q^2)$, we
obtain an infinite product in \eqref{rkibilat} when $i=k$:
$$
R_{k,k}(1,1/q;1;q^2) =
\frac{(-q)_{\infty}(q^{2k-1};q^{2k-1})_{\infty}}{(q)_{\infty}(-q^{2k-1};q^{2k-1})_{\infty}}.
$$
In the $r_{k,k}(s,t,m,n)$ of Theorem \ref{thm1}, the notion of an
unattached part then transfers to overpartitions by saying that an
odd part $2j-1$ occurs unattached if $2j, \overline{2j}$, and
$\overline{2j-1}$ do not occur.  The valuation
$v_{j}((\lambda,\mu))$ becomes a valuation defined on
overpartitions at even numbers, $v^1_{2j}(\lambda) =
f_{2j}(\lambda) + f_{\overline{2j-1}}(\lambda) +
f_{\overline{2j}}(\lambda) + \chi(\text{$2j-1$ occurs unattached
in $\lambda$})$.   The corresponding theorem is the overpartition
analogue of the Andrews-Gordon-G\"ollnitz identities mentioned in
the introduction:

\begin{corollary}[Lovejoy, \cite{Lo2}] \label{cor1}
Let $A^1_k(n)$ denote the number of overpartitions of $n$ into
parts not divisible by $2k-1$.  Let $B^1_k(n)$ denote the number
of overpartitions $\lambda$ of $n$ such that $(i)$ $v^1_2(\lambda)
\leq i-1$ and $(ii)$ for all $j \geq 1$ we have $f_{2j}(\lambda) +
v^1_{2j+2}(\lambda) \leq k-1$.  Then $A^1_k(n) = B^1_k(n)$.
\end{corollary}

All of the other results mentioned in the introduction (under
Theorem \ref{thm2}) follow in the same way.  Gordon's
generalization of the Rogers-Ramanujan identities \cite{Go1}
corresponds to the case $R_{k,i}(0,0;1;q)$,  Bressoud's
Rogers-Ramanujan identities for even moduli \cite{Br1} is the case
$\tilde{R}_{k,i}(0,0;1;q)$, Andrews' generalization of the
G\"ollnitz-Gordon identities \cite{An.1} is the case
$R_{k,i}(0,1/q;1;q^2)$, and the Gordon's theorems for
overpartitions \cite{Lo1} are the cases $R_{k,k}(0,1;1;q)$ and
$R_{k,1}(0,1/q;1;q)$.  The reader may work out the details, or
consult \cite{Co-Lo-Ma1,Co-Ma1}, where these are discussed along
with several other families of identities coming from the case
$a=0$.

As another example when neither $a$ nor $b$ is $0$, let us take $a
= \sqrt{-1}$ and $b = -\sqrt{-1}$ in \eqref{rkitildebilat}.  We
obtain an infinite product when $i = k-1$,
\begin{equation}
\tilde{R}_{k,k-1}(\sqrt{-1},-\sqrt{-1};1;q)
=
\frac{(-q)_{\infty}(-q^2;q^2)_{\infty}(q^{k-1};q^{k-1})_{\infty}}{(q)_{\infty}(q^2;q^2)_{\infty}(-q^{k-1};q^{k-1})_{\infty}}.
\end{equation}
Applying Theorem \ref{thm2} and splitting the generating functions
into real and imaginary parts, we obtain two weighted identities,
one of which is the following:
\begin{corollary} \label{cor2}
Let $A^2_k(n)$ denote the number of overpartition pairs where the
parts of $\mu$ are even and the parts of $\lambda$ are not
divisible by $k-1$.  Let $B^2_k(n)$ denote the number of
overpartition pairs of $n$ satisfying the conditions of Theorem
\ref{thm2} for $i = k-1$ (i.e., conditions $(i)$ and $(ii)$ of
Theorem \ref{thm1} and \eqref{three} if there is equality in
condition $(ii)$), having an even number of overlined parts, and
weighted by
\begin{equation} \label{weight}
(\sqrt{-1})^{\text{$\#$ overlined parts in
$\lambda$}}(-\sqrt{-1})^{\text{$\#$ overlined parts in $\mu$}}.
\end{equation}
Then $A^2_k(n) = B^2_k(n)$.
\end{corollary}
It is interesting to note that $A^2_k(n)$ is a non-weighted
counting function, while $B^2_k(n)$ is weighted.  Other
non-weighted interpretations of $A^2_k(n)$ may be found in
\cite{Lo3}.

For our last example, we consider the case $abq = 1$. The apparent
problem is that in this case there may be an unlimited number of
non-overlined zeros in $\mu$.  Indeed, the term
$1/(abxq)_{\infty}$ tends to infinity when $x =1$.  To remedy
this, we shall not consider $R_{k,i}(a,b;1;q)$ or
$\tilde{R}_{k,i}(a,b;1;q)$ as before, but the limits
\begin{equation} \label{limxto1}
\lim_{x \to 1} (1-x) R_{k,i}(a,b;x;q)
\end{equation}
and
\begin{equation} \label{limxto1tilde}
\lim_{x \to 1} (1-x) \tilde{R}_{k,i}(a,b;x;q).
\end{equation}
These limits follow easily from \eqref{jtp}:
\begin{equation} \label{limprod}
\lim_{x \to 1} (1-x) R_{k,i}(a,1/aq;x;q) =
\frac{(-aq,-1/a)_{\infty}(q^{i-1},q^{2k-i},q^{2k-1};q^{2k-1})_{\infty}}{(q)_{\infty}^2}.
\end{equation}
\begin{equation} \label{limprodtilde}
\lim_{x \to 1} (1-x) \tilde{R}_{k,i}(a,1/aq;x;q) =
\frac{(-aq,-1/a)_{\infty}(q^{i-1},q^{2k-i},q^{2k-1};q^{2k-1})_{\infty}}{(q)_{\infty}^2}.
\end{equation}
On the other hand, recalling Abel's lemma, which states that
$$
\lim_{x \to 1} (1-x)\sum_{n \geq 0}A_nx^n = \lim_{n \to \infty}
A_n,
$$
the limits \eqref{limprod} and \eqref{limprodtilde} may be
interpreted as the generating functions for those overpartition
pairs counted by $r_{k,i}(s,t,\infty,n)$ and
$\tilde{r}_{k,i}(s,t,\infty,n)$, respectively, the infinitude of
the number of parts corresponding to an infinite number of
non-overlined zeros in $\mu$.

For example, if we take $(a,b,q) = (1/q,1/q,q^2)$, then the
infinite product in \eqref{limprod} is
$$
\frac{(-q;q^2)_{\infty}^2(q^{2i-2},q^{4k-2i},q^{4k-2};q^{4k-2})_{\infty}}{(q^2;q^2)_{\infty}^2}.
$$
For the overpartition pairs of Theorem \ref{thm1}, those parts $j$
in $\lambda$ become $2j$, those parts $\overline{j}$ in $\lambda$
or $\mu$ become $2j-1$, and those parts $j$ in $\mu$ become
$2j-2$.  Given that the overlined parts are necessarily odd, the
overlining becomes redundant and we can talk simply about
partition pairs without repeated odd parts.  The definition of
unattached changes to: an even part $2j$ of $\mu$ is said to be
unattached in the partition pair $(\lambda,\mu)$ if $2j+1$ doesn't
occur in $\lambda$ or $\mu$ and $2j+2$ does not occur in $\mu$.
The valuation function becomes a valuation at even numbers,
$$
v^3_{2j}((\lambda,\mu)) = f_{2j}(\lambda) + f_{2j-1}(\lambda) +
f_{2j-1}(\mu) + \chi(\text{$2j-2$ occurs unattached in $\mu$}).
$$
We may then state:
\begin{corollary} \label{cor4}
For $i \geq 2$ let $A^3_{k,i}(n)$ denote the number of partition
pairs $(\lambda,\mu)$ of $n$ such that the odd parts cannot be
repeated, and the even parts of $\mu$ are not congruent to $0$ or
$\pm (2i-2)$ modulo $4k-2$. Let $B^3_{k,i}(n)$ denote the number
of partition pairs $(\lambda,\mu)$ of $n$ without repeated odd
parts, where $(i)$ $f_1(\lambda) + f_{2}(\lambda) + f_{1}(\mu)
\leq i-1$ and $(ii)$ for each $j \geq 1$ we have $f_{2j}(\lambda)
+ v^3_{2j+2}((\lambda,\mu)) \leq k-1$.  Then $A^3_{k,i}(n) =
B^3_{k,i}(n)$.
\end{corollary}
Of course, a similar result holds for $(a,b,q) = (1/q,1/q,q^2)$
using \eqref{limprodtilde} and Theorem \ref{thm2}.  This is left
to the interested reader.

\section{Lattice paths}
In the next three sections, we shall prove Theorems \ref{thm3} and
\ref{thm4}.  We begin in this section by defining the lattice
paths counted by $E_{k,i}(s,t,n)$ and $\tilde{E}_{k,i}(s,t,n)$ and
showing that their generating functions are \eqref{rkibilat} and
\eqref{rkitildebilat}, respectively.

We study paths in the first quadrant that use five kinds of
unitary steps,
\begin{itemize}
\item North-East (NE): $(x,y) \rightarrow (x+1,y+1)$, \item
South-East (SE): $(x,y) \rightarrow (x+1,y-1)$, \item South (S):
$(x,y) \rightarrow (x,y-1)$, \item South-West (SW): $(x,y)
\rightarrow (x-1,y-1)$, and \item East (E): $(x,0) \rightarrow
(x+1,0)$,
\end{itemize}
with the following additional restrictions:
\begin{itemize}
\item A South or South-West step can only appear after a
North-East step, \item An East step can only appear at height $0$,
\item The paths must start on the $y$-axis and end on the
$x$-axis.
\end{itemize}

The most important notion associated with these paths is the peak.
A \emph{peak} is a vertex preceded by a North-East step and
followed by a South step (in which case it can be labelled by $a$
or $b$ and called an $a$-peak or a $b$-peak, respectively), by a
South-West step (in which case it is called an $ab$-peak) or by a
South-East step (in which case it will be called a 1-peak). We say that a peak is \emph{marked by $a$} if it is an
$a$-peak or an $ab$-peak, and that it is \emph{marked by $b$} if
it is a $b$-peak or an $ab$-peak.

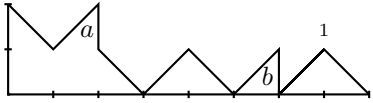
\begin{figure}[ht]
\psset{unit=0.6}
\begin{pspicture}(0,0)(8,2)
\psaxes[labels=none,ticksize=1pt](8,2)
\psline(0,2)(1,1)(2,2)(2,1)(3,0)(4,1)(5,0)(6,1)(6,0)(7,1)(6,0)(7,1)(8,0)
\uput[u](7,1){\tiny 1} \rput(1.75,1.4){\small $a$}
\rput(5.75,0.4){\small $b$}
\end{pspicture}
\caption{Example of a path. There are two $1$-peaks (located at
$(4,1)$ and $(7,1)$), an $a$-peak (located at $(2,2)$), a $b$-peak
(located at $(6,1)$), and an $ab$-peak (located at $(7,1)$). The
sequence of steps is SE-NE-S-SE-NE-SE-NE-S-NE-SW-NE-SE. The major
index is $2+4+6+7+7=26$.} \label{fig:expath1}
\end{figure}

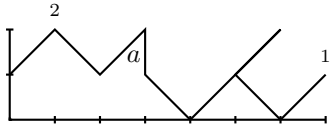
\begin{figure}[ht]
\psset{unit=0.6}
\begin{pspicture}(0,0)(7,2)
\psaxes[labels=none,ticksize=1pt](7,2)
\psline(0,1)(1,2)(2,1)(3,2)(3,1)(4,0)(6,2)(5,1)(6,0)(7,1)
\uput[u](1,2){\tiny 2} \uput[u](7,1){\tiny 1}
\rput(2.75,1.4){\small $a$}
\end{pspicture}
\caption{Another example. There is a $1$-peak (located at
$(1,2)$), an $a$-peak (located at $(3,2)$) and four $ab$-peaks
(located at $(1,2)$, $(1,2)$, $(6,2)$ and $(7,1)$). The sequence
of steps is NE-SW-NE-SW-NE-SE-NE-S-SE-NE-NE-SW-SE-NE-SW. The major
index is $1+1+1+3+6+7=19$.} \label{fig:expath2}
\end{figure}

The \emph{major index} of a path is the sum of the $x$-coordinates
of its peaks.  To avoid ambiguity in the graphical representation
of a path, we add a label to the $a$-peaks and $b$-peaks and we
may add a number above a vertex to indicate the presence of an
otherwise indistinguishable $ab$-peak (see Figures
\ref{fig:expath1} and \ref{fig:expath2} for examples).

We call these lattice paths \emph{generalized Bressoud-Burge
lattice paths}, for they are generalizations of some lattice paths
studied by Bressoud, \cite{Br3}, based on work of Burge
\cite{Bu1,Bu2}.  We might note that when there are no peaks marked
by $b$, we recover the paths studied in \cite{Co-Lo-Ma1} and
\cite{Co-Ma1}.

We now define the $(k,i)$-conditions that appear in Theorems
\ref{thm3} and \ref{thm4}.
\begin{definition}
We say that a path satisfies the \emph{odd $(k,i)$-conditions} if
it starts at height $k-i$ and its height is always less than $k$.
We say that a path satisfies the \emph{even $(k,i)$-conditions} if
it satisfies the odd $(k,i)$-conditions and for each peak of
coordinates $(x,k-1)$, we have $x-u+v \equiv i-1 \pmod 2$ where
$u$ is the number of $a$-peaks to the left of the peak and $v$ is
the number of $b$-peaks to the left of the peak.
\end{definition}

We first consider the paths satisfying the odd $(k,i)$-conditions.
Let $E_{k,i}(n,s,t,N)$ be the number of paths of major index $n$
with $N$ peaks, $s$ of which are marked by $a$ and $t$ of which
are marked by $b$, which satisfy the odd $(k,i)$-conditions. For
$0 < i \leq k$, let $\mathcal{E}_{k,i}(N)$ be the generating
function for these paths, that is
$\mathcal{E}_{k,i}(N)=\mathcal{E}_{k,i}(N,a,b,q)=\sum_{s,t,n}
E_{k,i}(n,s,t,N) a^s b^t q^n$. Moreover, for $0 \le i <k$,  let
$\Gamma_{k,i}(N)$ be the generating function for the paths
obtained by deleting the first NE step of a path which is counted
by $\mathcal{E}_{k,i+1}(N)$ and begins with a NE step.

\begin{proposition}\label{prop:recodd}
We have the following:
\begin{eqnarray}
\mathcal{E}_{k,i}(N) &=& q^N \Gamma_{k,i-1}(N) + q^N \mathcal{E}_{k,i+1}(N)  \quad\quad\text{$(0<i<k)$}\label{eq:recodd1},\\
\Gamma_{k,i}(N) &=& q^N \Gamma_{k,i-1}(N) + (a+b+q^{N-1}+abq^{1-N}) \mathcal{E}_{k,i+1}(N-1)
\quad\quad\text{$(0<i<k)$}\label{eq:recodd2},\\
\mathcal{E}_{k,k}(N) &=& q^N\Gamma_{k,k-1}(N) + q^N\mathcal{E}_{k,k}(N) \label{eq:recodd3},\\
\mathcal{E}_{k,i}(0) &=& 1\label{eq:recodd4},\\
\Gamma_{k,0}(N) &=& 0 \label{eq:recodd5}.
\end{eqnarray}
\end{proposition}
\begin{proof} The defining conditions of the generalized Bressoud-Burge paths imply that
the path has no peaks if and only if $N = 0$.  Hence
$\mathcal{E}_{k,i}(0) = 1$ (corresponding to the path that starts
at $(0,k-i)$ and descends with SE steps to $(k-i,0)$).  This is
\eqref{eq:recodd4}. If the path has at least one peak, then we take off
its first step and shift the path one unit to the left. If
$0<i<k$, then a path counted by $\mathcal{E}_{k,i}(N)$ starts with
a North-East step (corresponding to $q^N \Gamma_{k,i-1}(N)$) or a
South-East step (corresponding to $q^N \mathcal{E}_{k,i+1}(N)$).
This gives \eqref{eq:recodd1}.  For \eqref{eq:recodd2},
$\Gamma_{k,i}(N)$ is the generating function for the paths counted
by $\mathcal{E}_{k,i+1}(N)$ where the first North-East step was
deleted. These paths can start with a North-East step ($q^N
\Gamma_{k,i-1}(N)$), a South step ($(a+b)
\mathcal{E}_{k,i+1}(N-1)$), a South-East step ($q^{N-1}
\mathcal{E}_{k,i+1}(N-1)$) or a South-West step
($abq^{1-N}\mathcal{E}_{k,i+1}(N-1)$). If $i=k$ then a path
counted by $\mathcal{E}_{k,k}(N)$ starts with a North-East ($q^N
\Gamma_{k,k-1}(N)$) or an East step ($q^N \mathcal{E}_{k,k}(N)$).
The height of the paths is less than $k$, therefore no path which
starts at height $k-1$ can start with  a North-East step and so
$\Gamma_{k,0}(N)=0$.
\end{proof}

Notice that the recurrences and initial conditions above uniquely
define the generating functions $\mathcal{E}_{k,i}(N)$ and
$\Gamma_{k,i}(N)$. We shall exploit this fact to prove the
following generating functions:
\begin{theorem}\label{thm:gfNodd}
\begin{eqnarray*}
\mathcal{E}_{k,i}(N) &=& (ab)^N (-1/a,-1/b)_N q^N \sum_{n=-N}^N
\frac{(-1)^nq^{n((2k-1)n+3)/2+(k-i-1)n}}{(q)_{N-n}(q)_{N+n}}\\
\Gamma_{k,i}(N) &=& (ab)^N (-1/a,-1/b)_N \sum_{n=-N}^{N-1}
\frac{(-1)^nq^{n((2k-1)n+3)/2+(k-i-2)n}}{(q)_{N-n-1}(q)_{N+n}}
\end{eqnarray*}
\end{theorem}

\begin{proof}
Let
\begin{eqnarray*}
\mathcal E'_{k,i}(N)&=& (ab)^N (-1/a,-1/b)_N q^N \sum_{n=-N}^N
\frac{(-1)^nq^{n((2k-1)n+3)/2+(k-i-1)n}}{(q)_{N-n}(q)_{N+n}} \\
\Gamma'_{k,i}(N)&=& (ab)^N (-1/a,-1/b)_N \sum_{n=-N}^{N-1}
\frac{(-1)^nq^{n((2k-1)n+3)/2+(k-i-2)n}}{(q)_{N-n-1}(q)_{N+n}}.
\end{eqnarray*}
We will prove that these functions satisfy the five equations in
Proposition \ref{prop:recodd}.   We begin with \eqref{eq:recodd1}:
\begin{eqnarray*}
&& q^N\mathcal E'_{k,i+1}(N) + q^N \Gamma'_{k,i-1}(N)\\
&=& (ab)^N (-1/a,-1/b)_N q^N \sum_{n=-N}^N \frac{(-1)^nq^{n((2k-1)n+3)/2+(k-i-2)n}}{(q)_{N-n}(q)_{N+n}}q^N\\
&+& (ab)^N (-1/a,-1/b)_N \sum_{n=-N}^{N-1} \frac{(-1)^nq^{n((2k-1)n+3)/2+(k-i-1)n}}{(q)_{N-n-1}(q)_{N+n}}q^N\\
&=& (ab)^N (-1/a,-1/b)_N q^N \sum_{n=-N}^{N-1} \frac{(-1)^nq^{n((2k-1)n+3)/2+(k-i-1)n}}{(q)_{N-n}(q)_{N+n}}
\left(q^{N-n} + (1-q^{N-n}) \right)\\
&+& (ab)^N (-1/a,-1/b)^N q^N \frac{(-1)^Nq^{N((2k-1)N+3)/2+(k-i-1)N}}{(q)_0(q)_{2N}}\\
&=& (ab)^N (-1/a,-1/b)_N q^N \sum_{n=-N}^{N} \frac{(-1)^nq^{n((2k-1)n+3)/2+(k-i-1)n}}{(q)_{N-n}(q)_{N+n}}\\
&=& \mathcal E'_{k,i}(N).
\end{eqnarray*}
This gives \eqref{eq:recodd1}.  Notice that the above string of
equations holds for $i=k$. Next, we establish \eqref{eq:recodd2}:
\begin{eqnarray*}
&& q^N\Gamma'_{k,i-1}(N) + (a+b+q^{N-1}+q^{1-N}ab) \mathcal E'_{k,i+1}(N-1)\\
&=& (ab)^N (-1/a,-1/b)_N \sum_{n=-N}^{N-1}  \frac{(-1)^nq^{n((2k-1)n+3)/2+(k-i-1)n}}{(q)_{N-n-1}(q)_{N+n}}q^N\\
&+& (ab)^{N-1} (-1/a,-1/b)_{N-1} \sum_{n=-N+1}^{N-1} \frac{(-1)^nq^{n((2k-1)n+3)/2+(k-i-2)n}}{(q)_{N-n-1}(q)_{N+n-1}}
(a+b+q^{N-1}+q^{1-N}ab)\\
&=& (ab)^N (-1/a,-1/b)_N \sum_{n=-N}^{N-1} \frac{(-1)^nq^{n((2k-1)n+3)/2+(k-i-2)n}}{(q)_{N-n-1}(q)_{N+n}}q^{N+n}\\
&+& (ab)^{N-1} (-1/a,-1/b)_{N-1} \sum_{n=-N+1}^{N-1}  \frac{(-1)^nq^{n((2k-1)n+3)/2+(k-i-2)n}}{(q)_{N-n-1}(q)_{N+n-1}}\\
&&\quad\quad\quad\quad \times abq^{1-N}(1+a^{-1}q^{N-1})(1+b^{-1}q^{N-1})\\
&=& (ab)^N (-1/a,-1/b)_N \sum_{n=-N+1}^{N-1} \frac{(-1)^nq^{n((2k-1)n+3)/2+(k-i-2)n}}{(q)_{N-n-1}(q)_{N+n}}
\left( q^{N+n} + (1-q^{N+n}) \right)\\
&+& (ab)^N (-1/a,-1/b)_N \frac{(-1)^{-N}q^{-N((2k-1)(-N)+3)/2+(k-i-2)(-N)}}{(q)_{2N-1}(q)_0}\\
&=& (ab)^N (-1/a,-1/b)_N \sum_{n=-N}^{N-1} \frac{(-1)^nq^{n((2k-1)n+3)/2+(k-i-2)n}}{(q)_{N-n-1}(q)_{N+n}}\\
&=& \Gamma'_{k,i}(N).
\end{eqnarray*}
For \eqref{eq:recodd3}, we prove that $\mathcal E'_{k,k+1}(N) =
\mathcal E'_{k,k}(N)$ and then combine this with the fact that the
$\mathcal{E}'_{k,i}(N)$ satisfy \eqref{eq:recodd1} for $i=k$.
\begin{eqnarray*}
\mathcal E'_{k,k+1}(N) &=& (ab)^N (-1/a,-1/b)_N q^N \sum_{n=-N}^N \frac{(-1)^nq^{n((2k-1)n+3)/2-2n}}{(q)_{N-n}(q)_{N+n}}\\
&=& (ab)^N (-1/a,-1/b)_N q^N \sum_{n=-N}^N \frac{(-1)^nq^{-n((2k-1)(-n)+3)/2+2n}}{(q)_{N+n}(q)_{N-n}} \text{\hskip.2in
(replacing $n$ by $-n$)}\\
&=& (ab)^N (-1/a,-1/b)_N q^N \sum_{n=-N}^N \frac{(-1)^nq^{n((2k-1)n+3)/2-n}}{(q)_{N+n}(q)_{N-n}}\\
&=& \mathcal E'_{k,k}(N)
\end{eqnarray*}
Hence we have, using \eqref{eq:recodd1}:
$$
\mathcal E'_{k,k}(N) = q^N \mathcal E'_{k,k}(N) + q^N
\Gamma'_{k,k-1}(N).
$$
\par Notice that \eqref{eq:recodd4} is immediate.  Finally, for
\eqref{eq:recodd5}, we have
\begin{eqnarray*}
&&\Gamma'_{k,0}(N)\\
&=& (ab)^N (-1/a,-1/b)_N \sum_{n=-N}^{N-1} \frac{(-1)^nq^{n((2k-1)n+3)/2 +(k-2)n}}{(q)_{N-n-1}(q)_{N+n}}\\
&=& (ab)^N (-1/a,-1/b)_N \left( \sum_{n=0}^{N-1}
\frac{(-1)^nq^{n((2k-1)n+3)/2 +(k-2)n}}{(q)_{N-n-1}(q)_{N+n}} +
\sum_{n=-N}^{-1}  \frac{(-1)^nq^{n((2k-1)n+3)/2
+(k-2)n}}{(q)_{N-n-1}(q)_{N+n}} \right).
\end{eqnarray*}
Replacing $n$ by $-n-1$ in the second sum and simplifying gives
the negative of the first sum, which shows that $\Gamma'_{k,0}(N)
= 0$.  Now, since $\mathcal E'_{k,i}(N)$ and $\Gamma'_{k,i}(N)$
satisfy the same defining recurrences and initial conditions as
$\mathcal{E}_{k,i}(N)$ and $\Gamma_{k,i}(N)$, we have $\mathcal
E_{k,i}(N)=\mathcal E'_{k,i}(N)$ and
$\Gamma_{k,i}(N)=\Gamma'_{k,i}(N)$, which completes the proof.
\end{proof}


We are almost ready to prove that $E(s,t,n) = B(s,t,n)$ in Theorem
\ref{thm3}.  We just need a $q$-series lemma.
\begin{lemma}\label{prop:auxGauss}
For any integer $n$, we have
\begin{equation} \label{lemmaequation}
\sum_{N \ge |n|} \frac{(-q^n/a,-q^n/b)_{N-n} (abq)^{N-n}
(-aq,-bq)_n}{(q)_{N+n}(q)_{N-n}} =
\frac{(-aq,-bq)_\infty}{(q,abq)_\infty}
\end{equation}
\end{lemma}
\begin{proof}
We only prove the case $n \ge 0$. The case $n<0$ is identical, as
one may compute that
$$
\frac{(-q^{-n}/a,-q^{-n}/b)_{N+n} (abq)^{N+n}
(-aq,-bq)_{-n}}{(q)_{N-n}(q)_{N+n}} = \frac{(-q^n/a,-q^n/b)_{N-n}
(abq)^{N-n} (-aq,-bq)_n}{(q)_{N+n}(q)_{N-n}}.
$$
We have
\begin{eqnarray*}
&&\sum_{N \ge n} \frac{(-q^n/a,-q^n/b)_{N-n} (abq)^{N-n} (-aq,-bq)_n}{(q)_{N+n}(q)_{N-n}} \\
&=&\sum_{N \ge 0} \frac{(-q^n/a,-q^n/b)_N (abq)^N (-aq,-bq)_n}{(q)_{N+2n}(q)_N}\\
&=&\frac{(-aq,-bq)_n}{(q)_{2n}} \sum_{N \ge 0} \frac{(-q^n/a,-q^n/b)_N (abq)^N}{(q,q^{2n+1})_N}\\
&=&\frac{(-aq,-bq)_n}{(q)_{2n}} \frac{(-aq^{n+1},-bq^{n+1})_\infty}{(q^{2n+1},abq)_\infty}\\
&&\text{ by Corollary 2.4 of \cite{An6} with $n \to N$, $a \to -q^n/a$, $b \to -q^n/b$ and $c \to q^{2n+1}$}\\
&=&\frac{(-aq,-bq)_\infty}{(q,abq)_\infty}.
\end{eqnarray*}
\end{proof}

\emph{Proof of the case $B_{k,i}(s,t,n) = E_{k,i}(s,t,n)$ of
Theorem \ref{thm3}.}  Using the generating function from Theorem
\ref{thm:gfNodd} and summing on $N$ using Lemma
\ref{prop:auxGauss}, we have
\begin{eqnarray*}
\sum_{s,t,n \ge 0} E_{k,i}(s,t,n) a^s b^t q^n &=& \sum_{N \ge 0} \mathcal E_{k,i}(N)\\
&=& \sum_{N \ge 0} (ab)^N (-1/a,-1/b)_N q^N \sum_{n=-N}^N \frac{(-1)^nq^{n((2k-1)n+3)/2+(k-i-1)n}}{(q)_{N-n}(q)_{N+n}}\\
&=& \sum_{n=-\infty}^{\infty} (-1)^n q^{n((2k-1)n+3)/2+(k-i-1)n} \sum_{N \ge |n|} \frac{(ab)^N (-1/a,-1/b)_N q^N}
{(q)_{N-n}(q)_{N+n}}\\
&=& \sum_{n=-\infty}^{\infty} \frac{(-ab)^n(-1/a,-1/b)_n q^{n((2k-1)n+3)/2+(k-i)n}}{(-aq)_n (-bq)_n}\times\\
&&\quad\quad\quad\quad \sum_{N \ge |n|} \frac{(abq)^{N-n} (-q^n/a,-q^n/b)_{N-n}} {(q)_{N-n}(q)_{N+n}}\\
&=& \frac{(-aq)_\infty (-bq)_\infty}{(q)_\infty (abq)_\infty} \sum_{n=-\infty}^{\infty}
\frac{(-ab)^n(-1/a,-1/b)_n q^{n((2k-1)n+3)/2+(k-i)n}}{(-aq)_n (-bq)_n}\\
&=& R_{k,i}(a,b;1;q) \hskip.2in \text{(by \eqref{rkibilat})}\\
&=& \sum_{s,t,n \ge 0} B_{k,i}(s,t,n) a^s b^t q^n \hskip.2in
\text{(by Theorem \ref{thm1}).}
\end{eqnarray*}
Hence we have
$$
E_{k,i} (s,t,n) = B_{k,i} (s,t,n).
$$
\qed

We now treat the paths satisfying the even $(k,i)$-conditions.
Many of the arguments are similar to those for the paths
satisfying the odd $(k,i)$-conditions.  Hence, we shall not be as
verbose with details.  Let $\tilde{E}_{k,i}(n,s,t,N)$ be the
number of paths of major index $n$ with $N$ peaks, $s$ of which
are marked by $a$ and $t$ of which are marked by $b$, which
satisfy the even $(k,i)$-conditions. Let $\tilde{\mathcal
E}_{k,i}(N)$ and $\tilde{\Gamma}_{k,i}(N)$ be the even analogues
of $\mathcal{E}_{k,i}(N)$ and $\Gamma_{k,i}(N)$.

\begin{proposition}\label{prop:receven}
\begin{eqnarray}
\tilde{\mathcal{E}}_{k,i}(N) &=& q^N \tilde{\Gamma}_{k,i-1}(N) + q^N \tilde{\mathcal{E}}_{k,i+1}(N)   \quad\quad\text{$(0 < i <k)$},\\
\tilde{\Gamma}_{k,i}(N) &=& q^N \tilde{\Gamma}_{k,i-1}(N) +
(a+b+q^{N-1}+abq^{1-N}) \tilde{\mathcal{E}}_{k,i+1}(N-1)
\quad\quad\text{$(0 < i < k)$},\\
\tilde{\mathcal{E}}_{k,k}(N) &=& q^N \tilde{\mathcal E}_{k,k-1}(N) + q^N \tilde{\Gamma}_{k,k-1}(N),\\
\tilde{\mathcal{E}}_{k,i}(0) &=& 1,\\
\tilde{\Gamma}_{k,0}(N) &=& 0.
\end{eqnarray}
\end{proposition}
\begin{proof}
If $i < k$, we proceed just as in the proof of the Proposition
\ref{prop:recodd}.  If the path is not empty, then taking off its
first step increases or decreases $i$ by 1 and thus changes the
parity of $i-1$.  Moreover, all the peaks are shifted by 1, so the
parity of $x-u+v$ is not changed (for the recurrence for
$\tilde{\Gamma}_{k,i}(N)$, if the step we remove is a South step,
the peaks are not shifted but $u$ or $v$ decreases by 1 for all
peaks, so the result is the same).

The case $i=k$ needs further explanation.  The paths counted by
$\tilde{\mathcal{E}}_{k,k}(N)$ begin with either an East or a
North-East step.  Those that begin with a North-East step where
this step is deleted are the paths counted by
$\tilde{\Gamma}_{k,k-1}(N)$.  Shifting these one unit to the left
contributes the term $q^N\tilde{\Gamma}_{k,k-1}(N)$.

For the paths that begin with an East step, first observe that the
fact that every peak of coordinates $(x, k-1)$ satisfies $x-u+v
\equiv k-1 \pmod 2$ is equivalent to the fact that every peak of
coordinates $(x, k-1)$ has an even number of East steps to its
left.  We now consider two cases for the paths counted in
$\tilde{\mathcal E}_{k,k}(N)$ that start with an East step where
this step has been deleted.  If the path does not have any other
East step, then there is no peak of height $k-1$ and so we may
shift the path upward, i.e.\ each vertex of the path $(x, y)$ is
changed to $(x, y + 1)$.  Shifting to the left then creates a path
in $\tilde{\mathcal E}_{k,k-1}(N)$ that does not have any vertex
of the form $(x, 0)$. If the path does contain another East step,
then the path before the first of these other East steps is
shifted up, the East step is changed to a South-East step and the
rest of the path is not changed.  Shifting to the left creates a
path in $\tilde{\mathcal E}_{k,k-1}(N)$ that has at least one
vertex of the form $(x, 0)$.  This gives the term
$q^N\tilde{\mathcal{E}}_{k,k-1}(N)$.
\end{proof}

As in the odd case, the recurrences and initial conditions above
uniquely define the functions $\tilde{\mathcal{E}}_{k,i}(N)$ and
$\tilde{\Gamma}_{k,i}(N)$.  In this case, we have
\begin{theorem}\label{thm:gfNeven}
\begin{eqnarray*}
\tilde{\mathcal{E}}_{k,i}(N) &=& (ab)^N (-1/a,-1/b)_N q^N \sum_{n=-N}^N (-1)^n
\frac{q^{kn^2+(k-i-1)n - 2\binom{n}{2}}}{(q)_{N-n}(q)_{N+n}}\\
\tilde{\Gamma}_{k,i}(N) &=& (ab)^N (-1/a,-1/b)_N \sum_{n=-N}^{N-1} (-1)^n
\frac{q^{kn^2+(k-i-2)n - 2\binom{n}{2}}}{(q)_{N-n-1}(q)_{N+n}}
\end{eqnarray*}
\end{theorem}
The proof is omitted since it is very similar to that of Theorem \ref{thm:gfNodd}.

\emph{Proof of the case $\tilde{B}_{k,i}(s,t,n) =
\tilde{E}_{k,i}(s,t,n)$ of Theorem \ref{thm4}.}  This is identical
to the case of $B_{k,i}(s,t,n) = E_{k,i}(s,t,n)$ proven above.
Summing the generating function for $\tilde{\mathcal{E}}_{k,i}(N)$
over $N$ in Theorem \ref{thm:gfNeven}, changing the order of
summation and using Lemma \ref{prop:auxGauss} we get
\begin{eqnarray*}
\sum_{s,t,n \ge 0} \tilde{E}_{k,i}(s,t,n) a^s b^t q^n &=&
\tilde{R}_{k,i}(a,b;1;q) \\ &=& \sum_{s,t,n \ge 0}
\tilde{B}_{k,i}(s,t,n) a^s b^t q^n,
\end{eqnarray*}
and we conclude that
$$
\tilde{E}_{k,i} (s,t,n) = \tilde{B}_{k,i} (s,t,n).
$$
\qed

\section{Successive Ranks}
In this section we turn to the overpartition pairs counted by
$C_{k,i}(s,t,n)$ and $\tilde{C}_{k,i}(s,t,n)$.  We construct a
bijection between the relevant pairs and the lattice paths of the
previous section, which will establish the equality of
$C_{k,i}(s,t,n)$ and $E_{k,i}(s,t,n)$ (resp.
$\tilde{C}_{k,i}(s,t,n)$ and $\tilde{E}_{k,i}(s,t,n)$). This is a
generalization of overpartition-theoretic work in \cite{Co-Lo-Ma1}
and \cite{Co-Ma1}.

The Frobenius representation of an overpartition pair \cite{Co-Lo1,Co-Lo2,Lo3} of $n$ is a two-rowed array 
$$
\begin{pmatrix}
a_1 & a_2 & ... & a_N \\
b_1 & b_2 & ... & b_N
\end{pmatrix}
$$
where $(a_1,\ldots , a_N)$ and  $(b_1,\ldots ,b_N)$ are
overpartitions into nonnegative parts where $N + \sum (a_i + b_i)
= n$.  This is called the Frobenius representation of an
overpartition pair because these arrays are in bijection with
overpartition pairs of $n$ \cite{Co-Lo1,Ye1}.

We now define the successive ranks of an overpartition pair using
the Frobenius representation.
\begin{definition} \label{succrank}
If an overpartition pair has Frobenius representation
$$
\begin{pmatrix}
a_1&a_2&\cdots&a_N\\
b_1&b_2&\cdots&b_N\\
\end{pmatrix}
$$
then its $i$th successive rank $r_i$ is $a_i - b_i$ minus the number of non-overlined
parts in $\{b_{i+1},\ldots,b_N\}$ plus the number of non-overlined parts in $\{a_{i+1},\ldots,a_N\}$.
\end{definition}
For example, the successive ranks of
$\begin{pmatrix}
\overline{7} &4 &\overline{2} &0\\
\overline{3} &3 &1 &\overline{0}
\end{pmatrix}
$ are $(4,1,2,0)$.

We shall prove the following:
\begin{proposition}
There exists a one-to-one correspondence between the paths of major index $n$ counted by $E_{k,i}(s,t,n)$
and the overpartition pairs of $n$ counted by $C_{k,i}(s,t,n)$. This correspondence is such that the paths
have $N$ peaks if and only if the Frobenius representation of the overpartition pair has $N$ columns.
\label{prop:bijpathranksodd}
\end{proposition}
\begin{proof}
We prove this proposition by a direct mapping which is a
generalization of a mapping in \cite{Co-Ma1}. Given a lattice path
counted by $E_{k,i}(s,t,n)$, which starts at $(0,k-i)$, and a peak
$(x,y)$, let $u$ (resp. $v$) be the number of $a$-peaks to the
left of the peak (resp. the number of $b$-peaks to the left of the
peak). Starting on the left of the path, we construct a two-rowed
array from the right by mapping this peak to a column
$\begin{pmatrix} p
\\ q \end{pmatrix}$, where
\begin{equation} \label{evenp}
p = (x+k-i-y+u-v)/2
\end{equation}
and
\begin{equation} \label{evenq}
q = (x-k+i+y-2-u+v)/2,
\end{equation}
if there are an even number of East steps to the left of the peak, and
\begin{equation} \label{oddp}
p = (x+k-i+y-1+u-v)/2
\end{equation}
and
\begin{equation} \label{oddq}
q = (x-k+i-y-1-u+v)/2,
\end{equation}
if there are an odd number of East steps to the left of the peak.
Moreover, we overline the corresponding parts as follows:
\begin{itemize}
\item if the peak is a 1-peak, we overline $p$ and $q$,
\item if the peak is an $a$-peak, we overline $p$, and
\item if the peak is a $b$-peak, we overline $q$.
\end{itemize}
For example, the path counted by $E_{5,3}(3,4,115)$ in Figure
\ref{fig:expathranks} below maps to the overpartition pair counted
by $C_{5,3}(3,4,115)$ whose Frobenius representation is
$$
\begin{pmatrix}
14 &\overline{12} &12 &8 &\overline{7} &\overline{4} &\overline{3} &2\\
\overline{9} &\overline{8} &8 &\overline{7} &\overline{5}
&\overline{4} &3 &1
\end{pmatrix}.
$$

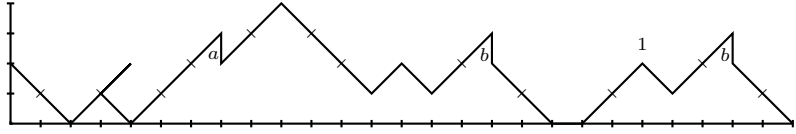
\begin{figure}[ht]
\psset{unit=0.4} \psaxes[labels=none,ticksize=1pt](26,4)
\begin{pspicture}(0,0)(26,4)
\psline(0,2)(2,0)(4,2)(3,1)(4,0)(7,3)(7,2)(9,4)(12,1)(13,2)(14,1)(16,3)(16,2)(18,0)(19,0)(21,2)(22,1)(24,3)(24,2)(26,0)
\psdots[dotstyle=x](1,1) \psdots[dotstyle=x](3,1)
\psdots[dotstyle=x](5,1) \psdots[dotstyle=x](6,2)
\psdots[dotstyle=x](8,3) \psdots[dotstyle=x](10,3)
\psdots[dotstyle=x](11,2) \psdots[dotstyle=x](15,2)
\psdots[dotstyle=x](17,1) \psdots[dotstyle=x](20,1)
\psdots[dotstyle=x](23,2) \psdots[dotstyle=x](25,1)
\uput[u](21,2){\tiny 1}
\rput(6.75,2.3){\tiny $a$} \rput(15.75,2.3){\tiny $b$}
\rput(23.75,2.3){\tiny $b$}
\end{pspicture}
\caption{A path counted by $E_{5,3}(3,4,115)$.}
\label{fig:expathranks}
\end{figure}

To establish the proposition, we must show that the result of this
mapping is indeed the Frobenius representation of an overpartition
pair counted by $C_{k,i}(s,t,n)$ and that the mapping is
invertible.  This is a somewhat tedious argument, involving a few
pages of small calculations and observations.  Ultimately we will
omit some details where these are similar to previous ones.

First, it may not even be clear that $p$ and $q$ defined above are
integers. To see this, note that at the starting point $(x,y) =
(0,k-i)$ of the path, the quantities $p$ and $q$ in \eqref{evenp}
and \eqref{evenq} are integers.  The parities of $x-y$, $x+y$,
$u-v$, and $u+v$ are preserved by $NE$, $SE$, and $SW$ steps,
while a $S$ step changes the parity of each of these. The only
problem is with an $E$ step, which changes the parity of $x-y$ and
$x+y$. This gives rise to the two cases for the definition in $p$
and $q$, which guarantees that the two-rowed array contains
integer entries.

Next, it is clear that the number of peaks in the path is equal to
the number of columns in the corresponding array.  It is also
clear that if the path contains $s$ (resp. $t$) peaks marked by
$a$ (resp. $b$), then the two-rowed array has $s$ (resp. $t$)
non-overlined parts in the bottom (resp. top) row.  Regarding $n$,
in either definition of $p$ and $q$ above, we have
\begin{equation} \label{pq1x}
p+q+1=x. \end{equation}
Hence, if $n$ is the major index of the
path, then $n$ is the sum of all entries of the corresponding
array and the number of columns.

Applying Definition \ref{succrank}, we compute the successive
ranks of the two-rowed array.  The peaks all have height at least
one, thus for a peak $(x,y)$ which is preceded by an even number
of East steps, we have:
\begin{eqnarray} \label{succrank1}
&&1 \le y=k-i+1+q-p+u-v\le k-1 \nonumber \\
&\Leftrightarrow&-i+2\le p-q-u+v \le k-i \\
&\Leftrightarrow&\text{the corresponding successive rank is $\ge
-i+2$ and $\le k-i$}, \nonumber
\end{eqnarray}
and if the peak is preceded by an odd number of East steps, we
have:
\begin{eqnarray} \label{succrank2}
&&1 \le y=p-q-u+v-k+i\le k-1 \nonumber \\
&\Leftrightarrow&k-i+1\le p-q-u+v \le 2k-i-1 \\
&\Leftrightarrow&\text{the corresponding successive rank is $\ge
k-i+1$ and $\le 2k-i-1$}. \nonumber
\end{eqnarray}
Hence, the successive ranks of the two-rowed array are all in the
interval $[-i+2,2k-i-1]$.

Finally, we need to prove that the two-rowed array we constructed
has an overpartition into non-negative parts in each row.  In what
follows, let $(x_j,y_j)$ be the coordinates of the $j$th peak from
the right and $\begin{pmatrix} p_j \\ q_j \end{pmatrix}$ be the
corresponding column, the $j$th column from the left.

First, we show that $p_N \ge 0$. If the leftmost peak has an even
number of East steps to its left, then $p_N = (x_N+k-i-y_N)/2$. It
is obvious that any vertex has a greater (or equal) value of $x-y$
than the previous vertex in the path. Since the path begins at
$(0,k-i)$, we have $x-y=-k+i$ at the beginning of the path and
thus we have $x-y \ge -k+i$ for all vertices and in particular for
the leftmost peak. Now if that peak has an odd number of East
steps to its left, then $p_N=(x_N+y_N+k-i-1)/2$. Since $x_N \ge 1$
and $y_N \ge 1$, we get that $p_N \ge 0$.

Next, we show that $q_N \ge 0$. This can be proven similarly.  If
the leftmost peak has an even number of East steps to the left,
then $q_N = (x_N + y_N - (k-i) - 2)/2$.  The path begins at
$(0,k-i)$, the only steps allowed before the first peak do not
decrease $x+y$, and there must be one $NE$ step before the first
peak, which increases $x+y$ by $2$.  Hence $x_N + y_N - (k-i) - 2
\geq 0$.  If the leftmost peak has an odd number of East steps to
the left, then $q_N = (x_N-y_N -(k-i) - 1)/2$.  Here the path
passes through the point $(k-i+1,0)$, and since $x-y$ never
decreases we have $q_N \geq 0$.

Having shown that all entries of the two-rowed array are
non-negative, we now argue that the sequences $\{p_j\}$ and
$\{q_j\}$ are overpartitions, i.e., that $p_j \ge p_{j+1}$ (resp.
$q_j \geq q_{j+1}$) with strict inequality if $p_{j+1}$ (resp.
$q_{j+1}$) is overlined. Let us show first that $p_j \ge p_{j+1}$.
We consider four cases. If the $j$th peak and the $j+1$th peak
both have an even number of East steps to their left, then $p_j -
p_{j+1} = (x_j - x_{j+1} - y_j + y_{j+1} + u_j - u_{j+1} - v_j +
v_{j+1})/2$. We always have $x_j -x_{j+1} \ge y_{j} - y_{j+1}$. We
can only have $u_j - u_{j+1} - v_j + v_{j+1} < 0$ if the $j+1$th
peak is a $b$-peak, but in that case we have $x_j -x_{j+1}
> y_{j} - y_{j+1}$. If the $j$th peak and the $j+1$th peak both
have an odd number of East steps to their left, the proof is
identical. If the $j$th peak has an odd number of East steps to
its left and the $j+1$th peak has an even number of East steps to
its left, the result is easily shown using the fact that $x_j -
x_{j+1} \ge 2$ since there is at least an East step between the
two peaks. In the final case, where the $j$th peak has an even
number of East steps to the left and the $j+1$th peak has an odd
number of East steps to its left, we have $p_j - p_{j+1} = (x_j -
x_{j+1} - y_j - y_{j+1} + 1 + u_j - u_{j+1} - v_j + v_{j+1})/2$.
Since there is at least one East step between the $j$th peak and
the $j+1$th peak, we have $x_j - x_{j+1} \ge y_j + y_{j+1}$ unless
the $j+1$th peak is an $ab$-peak (see Figure \ref{fig:aux01}).
Since $u_j - u_{j+1}$ can only be equal to 0 or 1 (the same holds
for $v_j - v_{j+1}$), we have $1 + u_j - u_{j+1} - v_j + v_{j+1}
\ge 0$ and therefore, $p_j - p_{j+1} \ge 0$.  If the $j+1$th peak
is an $ab$-peak, we have $x_j - x_{j+1} \ge y_j + y_{j+1} - 1$,
$u_j = u_{j+1}$ and $v_j = v_{j+1}$. Thus, we also have $p_j -
p_{j+1} \ge 0$.

\begin{figure}[ht]
\psset{unit=0.6}
\begin{pspicture}(0,-1)(9,4)
\psline(0,0)(9,0) \psline(0,3)(1,4)(1,3)(4,0) \psline(7,0)(9,2)
\psline[linestyle=dashed](1,3)(1,0)
\psline[linestyle=dashed](4,3)(4,0)
\psline[linestyle=dashed](7,2)(7,0)
\psline[linestyle=dashed](9,2)(9,0) \uput[d](1,0){$x_{j+1}$}
\uput[d](9,0){$x_j$} \psline[arrows=<->](1,3)(4,3)
\psline[arrows=<->](4,2)(7,2) \psline[arrows=<->](7,2)(9,2)
\uput[u](2.5,3){\small $\ge y_{j+1}-1$} \uput[u](5.5,2){$\ge 1$}
\uput[u](8,2){$\ge y_j$}
\end{pspicture}
\hspace{1cm}
\begin{pspicture}(0,-1)(9,5)
\psline(0,0)(9,0) \psline(1,5)(0,4)(4,0) \psline(7,0)(9,2)
\psline[linestyle=dashed](1,5)(1,0)
\psline[linestyle=dashed](4,3)(4,0)
\psline[linestyle=dashed](7,2)(7,0)
\psline[linestyle=dashed](9,2)(9,0) \uput[d](1,0){$x_{j+1}$}
\uput[d](9,0){$x_j$} \psline[arrows=<->](1,3)(4,3)
\psline[arrows=<->](4,2)(7,2) \psline[arrows=<->](7,2)(9,2)
\uput[u](2.5,3){\small $\ge y_{j+1}-2$} \uput[u](5.5,2){$\ge 1$}
\uput[u](8,2){$\ge y_j$}
\end{pspicture}
\caption{If the $j+1$th peak is not an $ab$-peak (left), we have
$x_j - x_{j+1} \ge y_j + y_{j+1}$. If the $j+1$th peak is an
$ab$-peak (right), we only have $x_j - x_{j+1} \ge y_j + y_{j+1} -
1$ but since $u_j = u_{j+1}$ and $v_j = v_{j+1}$, $p_j - p_{j+1}$
is indeed nonnegative.} \label{fig:aux01}
\end{figure}
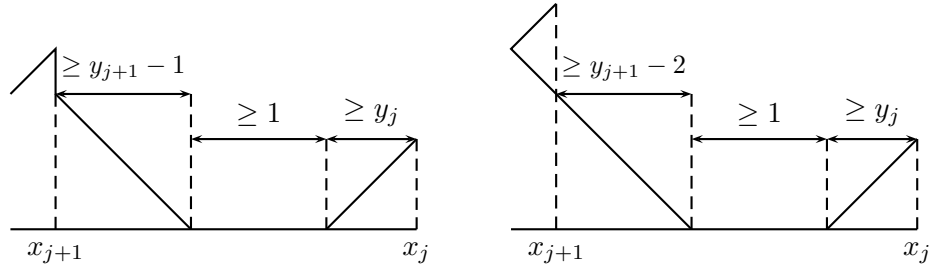

So, we have seen that in all cases we have $p_j \geq p_{j+1}$. If
$p_{j+1}$ is overlined, then the $j+1$th peak is a 1-peak or an
$a$-peak.  Going back through the above arguments, one finds that
the inequality $p_j \geq p_{j+1}$ is strict when the $j+1$th peak
is a $1$-peak or a $a$-peak.
The proof for the $\{q_j\}$ is quite similar to the case of the
$\{p_j\}$ above, so we omit the details.

That the mapping is invertible is rather straightforward.
Beginning at the left of the Frobenius representation of an
overpartition pair, at any column $\begin{pmatrix} p
\\ q \end{pmatrix}$ the values of $u$ and $v$ are determined by
the overlined parts in the columns to the right. Hence to recover
the location $(x,y)$ of a peak in the path from the column, we
need to solve either equations \eqref{evenp} and \eqref{evenq} or
equations \eqref{oddp} and \eqref{oddq} for $x$ and $y$.  Only one
set of these equations can be solved for positive $x$ and $y$.
Equation \eqref{pq1x} shows that $x$ does not increase as we
proceed. There is a unique way to fill in the steps between the
peaks, and the computations in \eqref{succrank1} and
\eqref{succrank2} show that the path never goes above height
$k-1$.  This completes the proof of the proposition.
\end{proof}

We will conclude this section by stating and proving the analogue
of Proposition \ref{prop:bijpathranksodd} for the functions
$\tilde{C}_{k,i}(s,t,n)$ and $\tilde{E}_{k,i}(s,t,n)$.

\begin{proposition}\label{prop:bijpathranskeven}
There exists a one-to-one correspondence between the paths of major index $n$ counted by $\tilde{E}_{k,i}(s,t,n)$
and the overpartition pairs of $n$ counted by $\tilde{C}_{k,i}(s,t,n)$. This correspondence is such that the paths
have $N$ peaks if and only if the Frobenius representation of the overpartition pair has $N$ columns.
\end{proposition}
\begin{proof}
A path counted by $\tilde{E}_{k,i}(s,t,n)$ is also counted by
$E_{k,i}(s,t,n)$ and an overpartition pair counted by
$\tilde{C}_{k,i}(s,t,n)$ is also counted by $C_{k,i}(s,t,n)$.
Hence we may apply the bijection used in the proof of Proposition
\ref{prop:bijpathranksodd} to a path counted by
$\tilde{E}_{k,i}(s,t,n)$.  We must then show that such paths
correspond to overpartition pairs where no successive rank can be
equal to $2k-i-1$. Indeed, if this was the case, we would have
$p-q-u+v=2k-i-1$ and from the map we know that $p-q-u+v=k-i-y+1$
or $k-i+y$. The first case is impossible when $k \geq 2$. The
second case implies that $y=k-1$ and $p = (x+u-v+2k-i-2)/2$. As
$p$ is an integer, we have $x-u+v \equiv i \pmod 2$. This is
forbidden by the last condition of the definition of
$\tilde{E}_{k,i}(s,t,n)$.
\end{proof}

\section{The Durfee dissection and a family of conjugations for overpartition pairs}
In this section we discuss the overpartition pairs counted by
$D_{k,i}(s,t,n)$ in Theorem \ref{thm3} and by
$\tilde{D}_{k,i}(s,t,n)$ in Theorem \ref{thm4}.  We complete our
proof of these two theorems using generating function identities
to show that these quantities are equal to $B_{k,i}(s,t,n)$ and
$\tilde{B}_{k,i}(s,t,n)$, respectively. The idea is to extend work
of Andrews \cite{An5} and Garvan \cite{Ga1} to overpartition pairs
via the Frobenius representation.

We begin by recalling a useful little bijection for
overpartitions, called the Joichi-Stanton algorithm \cite{Jo-St1}.
From an overpartition $\alpha$ into $N$ nonnegative parts, we
obtain a partition $\lambda$ into $N$ nonnegative parts and a
partition $\mu$ into distinct nonnegative parts less than $N$ as
follows: First, we initialize $\lambda$ to $\alpha$. Then, if the
$m$th part of $\alpha$ is overlined, we remove the overlining of
the $m$th part of $\lambda$, we decrease the $m-1$ first parts of
$\lambda$ by one and we add a part $m-1$ to $\mu$.

\begin{definition} \label{asso}
We say that $\lambda$ is the \emph{associated partition} of $\alpha$.
\end{definition}

Thus, given an overpartition pair we may decompose its Frobenius
representation into four partitions $\lambda_1$, $\mu_1$,
$\lambda_2$, $\mu_2$, where $\lambda_1$ and $\mu_1$ (resp.
$\lambda_2$ and $\mu_2$) are obtained by applying the
Joichi-Stanton algorithm to the top (resp. bottom) row.  For
example, the overpartition pair whose Frobenius representation is
$$
\pi =
\begin{pmatrix}
12 &12 &\overline{8} &7 &6 &\overline{3} &2 &\overline{1}\\
14 &12 &\overline{10} &\overline{8} &6 &5 &\overline{3} &2
\end{pmatrix}$$
gives $\lambda_1=(9,9,6,5,4,2,1,1)$, $\mu_1=(7,5,2)$,
$\lambda_2=(11,9,8,7,5,4,3,2)$ and $\mu_2=(6,3,2)$.

Next, we describe the notion of a $(k,i)$-admissible overpartition
pair occurring in the statement of Theorem \ref{thm3}.  This is
similar to, but not exactly the same as, the concept of
$(k,i)$-admissibility in \cite{An5}. Recall that the Durfee square
of a partition is the largest upper-left-justified square that
fits inside the Ferrers diagram of the partition \cite{An1}. Below
such a square, there is another partition and one may identify its
Durfee square, and so on, to obtain a sequence of successive
Durfee squares.

\begin{definition}
We say that an overpartition pair is $(k,i)$-admissible if the
conjugate, $\lambda '_2$, of the associated partition $\lambda_2$
of the bottom row of its Frobenius representation is obtained from
a partition $\nu$ into non-negative parts with at most $k-2$
Durfee squares by inserting a part of size $n_j$ into $\nu$ for
each $j$ with $i \leq j \leq k-1$. Here $n_j$ is the size of the
$j-1$th Durfee square of $\nu$, where the size of the $0$th Durfee
square is taken to be the number of columns in the Frobenius
representation of the overpartition pair.
\end{definition}

\begin{proposition}\label{prop:gfsquaresodd}
Recall the definition of $D_{k,i}(s,t,n)$ from Theorem \ref{thm3}.
We have the following generating function:
\begin{equation} \label{multiple1}
\sum_{s,t,n \geq 0} D_{k,i}(s,t,n)a^sb^tq^n =  \sum_{n_1 \ge
\cdots \ge n_{k-1} \ge 0} \frac{q^{n_1 + n_2^2 + \cdots +
n_{k-1}^2 + n_i + \cdots +
n_{k-1}}(-1/a,-1/b)_{n_{1}}a^{n_1}b^{n_1}}{(q)_{n_1-n_2}\cdots(q)_{n_{k-2}
- n_{k-1}}(q)_{n_{k-1}}}.
\end{equation}
\end{proposition}
\begin{proof}
Consider an overpartition pair counted by $D_{k,i}(s,t,n)$ whose
Frobenius representation has $n_1$ columns. By using the
Joichi-Stanton algorithm on each row, we can decompose our
overpartion pair in the following way:
\begin{itemize}
\item the top row, which is counted by
$$
\frac{(-1/b)_{n_1}b^{n_1}}{(q)_{n_1}},
$$
\item the partition $\mu_2$ into $n_1$ nonnegative parts coming
from the bottom row, which is counted by $(-1/a)_{n_1} a^{n_1}$,
\item the $n_1$ columns, which are counted by $q^{n_1}$, \item the
at most $k-2$ Durfee squares of a partition $\nu$, which are
counted by $q^{n_2^2 + \cdots +n_{k-1}^2}$, \item the regions
between the Durfee squares, which are counted by
$$
\qbinom{n_1}{n_2} \qbinom{n_2}{n_3} \cdots
\qbinom{n_{k-2}}{n_{k-1}},
$$
where
$$
\qbinom{n}{k}=\frac{(q)_n}{(q)_k(q)_{n-k}}
$$
is the generating function for partitions whose Ferrers diagrams
fit inside a $(n-k) \times k$ rectangle, and \item the inserted
parts, counted by
$$ q^{n_i + \cdots + n_{k-1}}.
$$
\end{itemize}
These last three together make up the conjugate $\lambda '_2$ of
the associated partition of the bottom row. Summing on $n_1,
\ldots, n_{k-1}$, we get the generating function:
\begin{eqnarray*}
&&\sum_{n_1 \ge n_2 \ge \cdots \ge n_{k-1} \ge 0}
\frac{(-1/b)_{n_1} b^{n_1}}{(q)_{n_1}} (-1/a)_{n_1} a^{n_1}
q^{n_1}
q^{n_2^2 + \cdots + n_{k-1}^2 + n_i + \cdots + n_{k-1}} \qbinom{n_1}{n_2} \cdots \qbinom{n_{k-2}}{n_{k-1}} \nonumber \\
&=& \sum_{n_1 \ge n_2 \ge \cdots \ge n_{k-1} \ge 0}
\frac{q^{n_1+n_2^2+\cdots+n_{k-1}^2 + n_i + \cdots +
n_{k-1}}(-1/a)_{n_1} a^{n_1} (-1/b)_{n_1}b^{n_1}} {(q)_{n_1-n_2}
\cdots (q)_{n_{k-2}-n_{k-1}} (q)_{n_{k-1}}}.
\end{eqnarray*}
\end{proof}

To incorporate $D_{k,i}(s,t,n)$ into Theorem \ref{thm3}, we use
the Bailey lattice structure from \cite{Ag-An-Br1} to transform
the generating function above to \eqref{rkibilat}.  Recall that a
pair of sequences $(\alpha_n, \beta_n)$ form a Bailey pair with
respect to $a$ if for all $n \geq 0$ we have
$$
\beta_n = \sum_{r = 0}^n \frac{\alpha_r}{(q)_{n-r}(aq)_{n+r}}.
$$
We shall employ the following lemma:
\begin{lemma} \label{lattice}
If $(\alpha_n,\beta_n)$ is a Bailey pair with respect to $q$, then
for all $0 \leq i \leq k$ we have
\begin{eqnarray} \label{latticeeq}
\frac{(abq)_{\infty}}{(q,-aq,-bq)_{\infty}} &\times& \sum_{n_1
\geq \cdots \geq n_k \geq 0} \frac{q^{n_1 + n_2^2 + \cdots n_k^2 +
n_{i+1} + \cdots
n_k}(-1/a,-1/b)_{n_1}(ab)^{n_1}}{(q)_{n_1-n_2}\cdots(q)_{n_{k-1} -
n_k}} \beta_{n_k} \nonumber \\
&=& \frac{\alpha_0}{(q)_{\infty}^2} +
\frac{1}{(q)_{\infty}^2}\sum_{n \geq 1}
\frac{(-1/a,-1/b)_n(ab)^nq^{(n^2-n)(i-1)+in}(1-q)}{(-aq,-bq)_n} \nonumber \\
&\times& \left(\frac{q^{(n^2+n)(k-i)}}{(1-q^{2n+1})}\alpha_n -
\frac{q^{((n-1)^2 + (n-1))(k-i)+ 2n-1}}{(1-q^{2n-1})} \alpha_{n-1}
\right)
\end{eqnarray}
\end{lemma}

\begin{proof}
This a special case of identity $(3.8)$ in \cite{Ag-An-Br1}.
Specifically, we let $a = q$, $\rho_1 = -1/a$, $\sigma_1 = -1/b$,
and then let $n$ as well as all remaining $\rho_i$ and $\sigma_i$
tend to $\infty$ in that identity to obtain \eqref{latticeeq}.
\end{proof}

\emph{Proof of the case $B_{k,i}(s,t,n) = D_{k,i}(s,t,n)$ of
Theorem \ref{thm3}.}  We use the Bailey pair with respect to $q$
\cite[p.468, (B3)]{Sl1},
$$
\beta_n = \frac{1}{(q)_{\infty}} \hskip.5in \text{and} \hskip.5in
\alpha_n = \frac{(-1)^nq^{n(3n+1)/2}(1-q^{2n+1})}{(1-q)}.
$$
Substituting into Lemma \ref{lattice} and simplifying, we obtain
\begin{align*}
\sum_{n_1 \ge \cdots \ge n_{k} \ge 0} & \frac{q^{n_1 + n_2^2 +
\cdots + n_{k}^2 + n_{i+1} + \cdots +
n_{k}}(-1/a,-1/b)_{n_{1}}a^{n_1}b^{n_1}}{(q)_{n_1-n_2}\cdots(q)_{n_{k-1}
- n_{k}} (q)_{n_k}} \\ =&
\frac{(-aq,-bq)_{\infty}}{(q,abq)_{\infty}} \left(1 + \sum_{n \geq
1} \frac{q^{kn^2+(k-i+1)n +
n(n+1)/2}(-ab)^n(-1/a,-1/b)_n}{(-aq,-bq)_n} \right. \\  +& \left.
\sum_{n \geq 1} \frac{q^{kn^2 -(k-i)n +
n(n+1)/2}(-ab)^n(-1/a,-1/b)_n}{(-aq,-bq)_n} \right).
\end{align*}
Replacing $n$ by $-n$ in the second sum, simplifying, and then
replacing $k$ by $k-1$ and $i$ by $i-1$ gives \eqref{rkibilat}.
\qed

Now we turn to the function $\tilde{D}_{k,i}(s,t,n)$.  We define
an operation on overpartition pairs, called $k$-conjugation, again
using the Frobenius representation. Recall the decomposition of
such a representation described after Definition \ref{asso}.  Let
$\lambda '_1$ (resp. $\lambda '_2$) be the conjugate of
$\lambda_1$ (resp. of $\lambda_2$).  Thus, $\lambda '_1$ and
$\lambda '_2$ are partitions into parts less than or equal to
$n_1$, where $n_1$ is the number of columns in the Frobenius
representation.  We consider two regions.  As above, we say that
the $0$th Durfee square of a partition has size $n_1$. The first
region $G_2$ we consider is the portion of $\lambda '_2$ below its
$(k-2)$-th Durfee square. The second region $G_1$ consists of the
parts of $\lambda '_1$ which are less than or equal to the size of
the $(k-2)$-th Durfee square of $\lambda '_2$.

\begin{definition}
For the \emph{$k$-conjugation} of an overpartition pair, we first
interchange these two regions $G_1$ and $G_2$ of $\lambda '_1$ and
$\lambda '_2$ to get two new partitions $\lambda ''_1$ and
$\lambda ''_2$.  Next, we conjugate these to get $\lambda '''_1$
and $\lambda '''_2$. Finally, we use the Joichi-Stanton algorithm
to assemble $\lambda '''_1$ and $\mu_1$ into the top row and
$\lambda '''_2$ and $\mu_2$ into the bottom row.
\end{definition}

We remark that if $\lambda '_2$ has less than $k-2$ Durfee
squares, the $k$-conjugation is the identity. Note that this
$k$-conjugation is a generalization of the $k$-conjugation for
overpartitions defined by Corteel and the present authors in
\cite{Co-Lo-Ma1} (which in turn was a generalization of Garvan's
$k$-conjugation for partitions \cite{Ga1}).

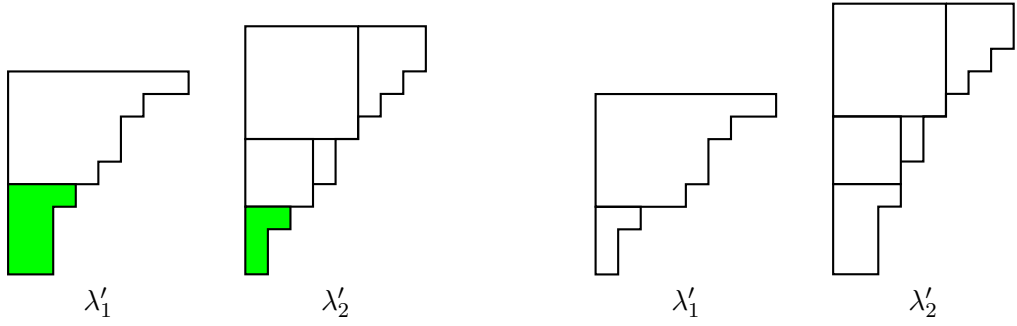
\begin{figure}[ht]
\psset{unit=0.3}
\begin{pspicture}(8,-11)
\pspolygon[dimen=middle](0,0)(8,0)(8,-1)(6,-1)(6,-2)(5,-2)(5,-4)(4,-4)(4,-5)(0,-5)
\pspolygon[dimen=middle,fillstyle=solid,fillcolor=green](0,-5)(3,-5)(3,-6)(2,-6)(2,-9)(0,-9)
\uput[d](4,-9){$\lambda'_1$}
\end{pspicture}
\hspace{0.5cm}
\begin{pspicture}(8,-13)
\psframe[dimen=middle](0,0)(5,-5)
\psframe[dimen=middle](0,-5)(3,-8)
\psline[dimen=middle](5,0)(8,0)(8,-2)(7,-2)(7,-3)(6,-3)(6,-4)(5,-4)(5,-5)(4,-5)(4,-7)(3,-7)
\pspolygon[dimen=middle,fillstyle=solid,fillcolor=green](0,-8)(2,-8)(2,-9)(1,-9)(1,-11)(0,-11)
\uput[d](4,-11){$\lambda'_2$}
\end{pspicture}
\hspace{2cm}
\begin{pspicture}(8,-10)
\pspolygon[dimen=middle](0,0)(8,0)(8,-1)(6,-1)(6,-2)(5,-2)(5,-4)(4,-4)(4,-5)(0,-5)
\pspolygon[dimen=middle](0,-5)(2,-5)(2,-6)(1,-6)(1,-8)(0,-8)
\uput[d](4,-8){$\lambda'_1$}
\end{pspicture}
\hspace{0.5cm}
\begin{pspicture}(8,-14)
\psframe[dimen=middle](0,0)(5,-5)
\psframe[dimen=middle](0,-5)(3,-8)
\psline[dimen=middle](5,0)(8,0)(8,-2)(7,-2)(7,-3)(6,-3)(6,-4)(5,-4)(5,-5)(4,-5)(4,-7)(3,-7)
\pspolygon[dimen=middle](0,-8)(3,-8)(3,-9)(2,-9)(2,-12)(0,-12)
\uput[d](4,-12){$\lambda'_2$}
\end{pspicture}
\caption{Illustration of the $4$-conjugation. For the initial
overpartition $\pi$, we have $\lambda '_1=(8,6,5,5,4,3,2,2,2)$ and
$\lambda '_2=(8,8,7,6,5,4,4,3,2,1,1)$. The regions highlighted are
interchanged by 4-conjugation, which gives $\lambda
'_1=(8,6,5,5,4,2,1,1)$ and $\lambda '_2=(8,8,7,6,5,4,4,3,3,2,2,2)$
for $\pi^{(4)}$, the $4$-conjugate of $\pi$.} \label{fig:kconj}
\end{figure}

Continuing with the example from after Definition \ref{asso}, it
is easy to see that we have $\lambda'_1=(8,6,5,5,4,3,2,2,2)$ and
$\lambda'_2=(8,8,7,6,5,4,4,3,2,1,1)$. For $k=4$, if we interchange
the two regions defined above, we get
$\lambda''_1=(8,6,5,5,4,2,1,1)$ and
$\lambda''_2=(8,8,7,6,5,4,4,3,3,2,2,2)$ (see Figure
\ref{fig:kconj}). Conjugating, we get
$\lambda'''_1=(8,6,5,5,4,2,1,1)$ and
$\lambda'''_2=(12,12,9,7,5,4,3,2)$. By applying the Joichi-Stanton
algorithm in reverse (remember that $\mu_1=(7,5,2)$ and
$\mu_2=(6,3,2)$), we see that the 4-conjugate of $\pi$ is
$$
\pi^{(4)} = \begin{pmatrix}
11 &9 &\overline{7} &7 &6 &\overline{3} &2 &\overline{1}\\
15 &15 &\overline{11} &\overline{8} &6 &5 &\overline{3} &2
\end{pmatrix}.
$$

\begin{definition}
We say that an overpartition pair is \emph{self-$k$-conjugate} if it is fixed by $k$-con\-ju\-ga\-tion.
\end{definition}

\begin{proposition}\label{prop:gfselfconj}
The generating function for self-$k$-conjugate overpartition pairs is
\begin{equation} \label{kselfconjgf}
\sum_{n_1 \ge n_2 \ge \cdots \ge n_{k-1} \ge 0}
\frac{q^{n_1+n_2^2+ \cdots+n_{k-1}^2}(-1/a)_{n_1}
a^{n_1}(-1/b)_{n_1}b^{n_1}}{(q)_{n_1-n_2} \cdots
(q)_{n_{k-2}-n_{k-1}} (q^2;q^2)_{n_{k-1}}},
\end{equation}
where $n_1$ is the number of columns of the Frobenius symbol and $n_2,\ldots,n_{k-1}$ are the sizes
of the $k-2$ first successive Durfee squares of $\lambda '_2$.
\end{proposition}
\begin{proof}
The decomposition of a self-$k$-conjugate overpartition pair is
similar to the decomposition of a $(k,k)$-admissible overpartition
pair.  We have the following pieces:
\begin{itemize}
\item $\mu_1$, which is counted by $(-1/b)_{n_1}b^{n_1}$, \item
$\mu_2$, which is counted by $(-1/a)_{n_1}a^{n_1}$, \item The
$n_1$ columns, which are counted by $q^{n_1}$, \item the $k-2$
Durfee squares of $\lambda '_2$, which are counted by
$q^{n_2^2+\cdots+n_{k-1}^2}$, \item the regions between the Durfee
squares of $\lambda '_2$, which are counted by
$$\qbinom{n_1}{n_2} \cdots \qbinom{n_{k-2}}{n_{k-1}},$$
\item the parts in $\lambda '_1$ which are $>n_{k-1}$ and of
course $\le n_1$: they are counted by
$$\frac{1}{(1-q^{n_{k-1}+1})\cdots(1-q^{n_1})} = \frac{(q)_{n_{k-1}}}{(q)_{n_1}},$$
\item the two identical regions $G_1$ and $G_2$, which are counted
by
$$\frac{1}{(q^2;q^2)_{n_{k-1}}}.$$
\end{itemize}
For example, in Figure \ref{fig:kconj} we do not have a
self-$4$-conjugate overpartition pair because the shaded regions
are not identical.

Summing on $n_1, n_2, \ldots, n_{k-1}$, we get the generating
function:
\begin{align*}
&\sum_{n_1 \ge n_2 \ge \cdots \ge n_{k-1} \ge 0}
(-1/b)_{n_1}b^{n_1} (-1/a)_{n_1} a^{n_1} q^{n_1}
q^{n_2^2+\cdots+n_{k-1}^2} \qbinom{n_1}{n_2} \cdots
\qbinom{n_{k-2}}{n_{k-1}} \frac{(q)_{n_{k-1}}}{(q)_{n_1}}
\frac{1}{(q^2;q^2)_{n_{k-1}}}\\
=& \sum_{n_1 \ge n_2 \ge \cdots \ge n_{k-1} \ge 0}
\frac{q^{n_1+n_2^2+\cdots+n_{k-1}^2}(-1/a)_{n_1} a^{n_1}
(-1/b)_{n_1}b^{n_1}} {(q)_{n_1-n_2} \cdots (q)_{n_{k-2}-n_{k-1}}
(q^2;q^2)_{n_{k-1}}}.
\end{align*}
\end{proof}


\begin{definition} \label{selfki}
We say that an overpartition pair is \emph{self-$(k,i)$-conjugate}
if it is obtained by taking a self-$k$-conjugate overpartition
pair and adding a part $n_j$ ($n_j$ is the size of the $(j-1)$-th
successive Durfee square of $\lambda '_2$) to $\lambda '_2$ for $i
\le j \le k-1$.
\end{definition}

Remember that we denote by $\tilde{D}_{k,i}(s,t,n)$ the number of
self-$(k,i)$-conjugate overpartition pairs of $n$ whose Frobenius
representations have $s$ non-overlined parts in their bottom rows
and $t$ non-overlined parts in their top rows.  We may now
complete the proof of Theorem \ref{thm4}.

\emph{Proof of the case $\tilde{D}_{k,i}(s,t,n) =
\tilde{B}_{k,i}(s,t,n)$ of Theorem \ref{thm4}.} It is obvious from
Proposition \ref{prop:gfselfconj} and Definition \ref{selfki} that
\begin{equation} \label{multiple2}
\sum_{s,t,n \geq 0} \tilde{D}_{k,i}(s,t,n) a^s b^t q^n = \sum_{n_1
\ge n_2 \ge \cdots \ge n_{k-1} \ge 0}
\frac{q^{n_1+n_2^2+\cdots+n_{k-1}^2+n_i+\cdots+n_{k-1}}(-1/a,-1/b)_{n_1}
a^{n_1}b^{n_1}}{(q)_{n_1-n_2} \cdots (q)_{n_{k-2}-n_{k-1}}
(q^2;q^2)_{n_{k-1}}}.
\end{equation}
Consider the Bailey pair with respect to $q$ \cite[p.468,
(E3)]{Sl1},
$$
\beta_n = \frac{1}{(q^2;q^2)_{\infty}} \hskip.5in \text{and}
\hskip.5in \alpha_n = \frac{(-1)^nq^{n^2}(1-q^{2n+1})}{(1-q)}.
$$
Substituting into Lemma \ref{lattice} and arguing in the case of
$D_{k,i}(s,t,n)$ above shows that \eqref{multiple2} is equal to
\eqref{rkitildebilat}.

\qed

\section{Concluding Remarks}
We wish to close with a look at some possible future research
topics.  First, several authors \cite{An.05,Br1,Lo3} have derived
combinatorial identities from Andrews' $J_{1,k,i}(a;x;q)$ when $k$
and/or $i$ are half-integers.  Can this idea be applied to the
$R_{k,i}(a,b;x;q)$ or $\tilde{R}_{k,i}(a,b;x;q)$?   For example,
we might mention that the $R_{2,k,3/2}(-q,-q^2;1;q^2)$ are
expressible as infinite products.

Second, it would be worthwhile to develop the recurrences for a
``tilde version" of Andrews' $J_{\lambda,k,i}$ for all $\lambda$
and see if there is perhaps an analogue of Andrews' general
Rogers-Ramanujan theorem \cite{An3}. Such a theorem was in fact
predicted by Bressoud \cite[p.19]{Br2}. Moreover, there are other
nice applications of Andrews' functions besides proving
combinatorial theorems. For instance, they have been used to prove
$q$-series identities of the Rogers-Ramanujan type \cite{An1,An4}
and in the study of $q$-continued fractions \cite{An1,An-Be1}. The
tilde analogues would, no doubt, be equally fruitful.

Finally, we now know that Andrews' functions $J_{0,k,i}(-;x;q)$
are generating functions for certain partitions, the
$J_{1,k,i}(-1/a;x;q)$ are generating functions for certain
overpartitions, and the $J_{2,k,i}(-1/a,-1/b;x;q)/(abxq)_{\infty}$
are generating functions for certain overpartition pairs.  The
natural question, of course, is what happens to all of the
combinatorial objects considered in this paper when we pass to
$\lambda = 3$?  There are a number of barriers that make it
unclear how to go beyond overpartition pairs.  First, from the
perspective of generating functions, passing from partitions to
overpartitions to overpartition pairs involves passing from
$q^{n^2}$ to $(-1/a)_na^nq^{n(n+1)/2}$ to $(-1/a,-1/b)_n(abq)^n$.
What would be next?  Second, in terms of Frobenius symbols, we
pass from symbols with partitions into distinct parts in both rows
to symbols with an overpartition in one row and a partition into
distinct parts in the other to symbols with overpartitions in both
rows.  Again, what would be next?  Third, in terms of the lattice
paths, a peak can be open when dealing with partitions, half-open
when we allow overpartitions, or closed when we pass to
overpartition pairs.  What would happen to these peaks and paths
in the next case?

In terms of $q$-series identities, we know that the correspondence
between an overpartition pair and its Frobenius symbol is the
essence of two famous identities, the $q$-Gauss summation and the
$_1 \psi_1$ summation \cite{Co1,Co-Lo1,Ye1}.  Perhaps a clue to
going beyond overpartition pairs lies in finding a natural
bijective proof of some generalization of the $q$-Gauss summation.
The $_6 \phi _5$ summation, for example, would be a good
candidate.  Some work toward a bijective proof of this identity is
presented in \cite{Co1}.

\end{document}